\documentclass[onefignum,onetabnum]{siamonline181217}

\usepackage{lipsum}
\usepackage{amsfonts}
\usepackage{graphicx}
\usepackage{epstopdf}
\ifpdf
  \DeclareGraphicsExtensions{.eps,.pdf,.png,.jpg}
\else
  \DeclareGraphicsExtensions{.eps}
\fi

\usepackage{enumitem}
\setlist[enumerate]{leftmargin=.5in}
\setlist[itemize]{leftmargin=.5in}

\newsiamremark{remark}{Remark}
\newsiamremark{example}{Example}
\newsiamthm{assumptions}{Assumptions}

\newsiamremark{hypothesis}{Hypothesis}
\crefname{hypothesis}{Hypothesis}{Hypotheses}
\crefname{remark}{Remark}{Remarks}
\crefname{example}{Example}{Examples}

\crefname{assumptions}{Assumptions}{Assumptions}

\newsiamthm{claim}{Claim}

\headers{On the well-posedness of Bayesian inverse problems}{J. Latz}

\usepackage{amsopn}

\ifpdf
\hypersetup{
  pdftitle={On the well-posedness of Bayesian inverse problems},
  pdfauthor={J. Latz}
}
\fi

\usepackage{latexsym,amsfonts,amsmath,theorem,amssymb}
\usepackage{pstricks,graphicx}
\usepackage{a4wide}
\usepackage{color}
\usepackage{caption}
\usepackage{theorem}
\usepackage{framed}
\usepackage{blindtext}
\usepackage{centernot}
\usepackage{todonotes}
\usepackage{rotating}
\usepackage{float}
\usepackage[caption=false]{subfig}
\usepackage{tikz}

\newcommand{\Prob}{\mathrm{Prob}}
\newcommand{\PP}{\mathrm{P}}

\newcommand{\dHel}{\mathrm{d}_{\mathrm{Hel}}}
\newcommand{\dtv}{\mathrm{d}_{\mathrm{tv}}}
\newcommand{\dProk}{\mathrm{d}_{\mathrm{Prok}}}
\newcommand{\dWas}{\mathrm{d}_{\mathrm{Was(p)}}}
\newcommand{\DKL}{\mathrm{D}_{\mathrm{KL}}}

\newcommand{\munoise}{\mu_{\mathrm{noise}}}
\newcommand{\mupost}{\mu_{\mathrm{post}}}
\newcommand{\pipost}{\pi_{\mathrm{post}}}
\newcommand{\muprio}{\mu_{\mathrm{prior}}}
\newcommand{\piprio}{\pi_{\mathrm{prior}}}

\newcommand{\vertiii}[1]{{\left\vert\kern-0.25ex\left\vert\kern-0.25ex\left\vert #1 
\right\vert\kern-0.25ex\right\vert\kern-0.25ex\right\vert}}

\graphicspath{{./figures/}}

\title{On the well-posedness of Bayesian inverse problems\thanks{Submitted to the editors DATE.
\funding{The author gratefully acknowledges the support by DFG and Technische Universit\"at M\"unchen through the International Graduate School of Science and Engineering within project 10.02 BAYES.}}}

\author{Jonas Latz\thanks{Zentrum Mathematik, Technische Universit\"{a}t M\"{u}nchen, Boltzmannstr. 3, 85747 Garching b.M., Germany
(jonas.latz@ma.tum.de).} 
}

\begin{document}
\maketitle

\begin{abstract}
The subject of this article is the introduction of a new concept of well-posedness of Bayesian inverse problems. 
The conventional concept of (Lipschitz, Hellinger) well-posedness in [Stuart 2010, Acta Numerica 19, pp. 451-559] is difficult to verify in practice and may be inappropriate in some contexts.
Our concept simply replaces the Lipschitz continuity of the posterior measure in the Hellinger distance by continuity in an appropriate distance between probability measures.
Aside from the Hellinger distance, we investigate well-posedness with respect to  weak convergence, the total variation distance, the Wasserstein distance, and also the Kullback--Leibler divergence.
We demonstrate that the weakening to continuity is tolerable and that the generalisation to other distances is important.
The main results of this article are proofs of well-posedness with respect to some of the aforementioned distances for large classes of Bayesian inverse problems. Here, little or no information about the underlying model is necessary; making these results particularly interesting for practitioners using black-box models.
We illustrate our findings with numerical examples motivated from machine learning and image processing.
\end{abstract}



\section{Introduction}

The representation of systems and processes in nature and technology with mathematical and computational models is fundamental in modern sciences and engineering.
For a  partially observable process, the \emph{model calibration} or \emph{inverse problem} is of particular interest.
It consists in fitting model parameters such that the model represents the underlying process.
Aside from classical mathematical models, such as partial differential equations or dynamical systems, inverse problems also play a central role in machine learning applications, for example, classification with deep neural networks or (non-)linear support vector machines, as well as regression with deep Gaussian processes.

The solvability of inverse problems is usually classified in terms of their \emph{well-posedness}.
In 1902, Hadamard \cite{Hadamard1902} defined:
\begin{equation*}
\text{\emph{``[\ldots] ces probl\'emes [\ldots] bien pos\'e, je veux dire comme possible et d\'etermin\'e."}}
\end{equation*}
According to Hadamard, a problem is \emph{well-posed} if the solution is \emph{possible} and  \emph{determined}; it can be \emph{found} and is \emph{exact}.
Today, we interpret this principle as: the solution of the inverse problem exists, is unique, and depends continuously on the data. 
The continuity in the data implies stability. 
The existence and stability allow us to find the solution (\emph{`possible'}) and uniqueness makes the solution exact (\emph{`d\'etermin\'e'}).
This is a justification for well-posedness from an analytical and computational viewpoint. 
From a statistical viewpoint, well-posedness does not only allow us to find the estimate; 
it also gives us robustness of the estimate with respect to marginal perturbations in the data:  
Since we often know that data is not precise, we should anticipate to see only marginal changes in the estimate with respect to marginal changes in the data. Otherwise, we should not consider the estimate trustworthy.

Measurement noise, complexity of the model, and a lack of data  typically lead to \emph{ill-posed} (i.e. not well-posed) inverse problems. 
The regularised least squares approach can sometimes be used to transform an ill-posed problem into a well-posed problem.
We refer to \cite{Chavent2010,Groetsch2015} for an introduction of that approach; however, it is outside of the scope of this article.

The Bayesian approach to inverse problems represents the uncertain model parameter as a random variable.
The random variable is distributed according to a \emph{prior (measure)}, which reflects uncertainty in the parameter.
Observing the data is then an event with respect to which the prior shall be conditioned.
The solution of the Bayesian inverse problem is this conditional probability measure, called the \emph{posterior (measure)}.
Stuart \cite{Stuart2010} transferred Hadamard's principle of   well-posedness to Bayesian inverse problems: 
the posterior exists, it is unique, and it is locally Lipschitz continuous with respect to the data.
In Stuart \cite{Stuart2010}, the distance between posteriors is measured in the Hellinger distance.
Several authors have discussed, what we choose to call, \emph{(Lipschitz, Hellinger) well-posedness} for a variety of Bayesian inverse problems. For example, elliptic partial differential equations \cite{Dashti2011,Iglesias2014}, level-set inversion \cite{ILS16}, Helmholtz source identification with Dirac sources \cite{Engel2018}, a Cahn-Hilliard model for tumour growth \cite{Kahle2019}, hierarchical prior measures \cite{LEU19}, 
stable priors in quasi-Banach spaces \cite{Sullivan2017_Wellposed,Sullivan2017}, convex and heavy-tailed priors \cite{Hosseini2017,Hosseini2017b}.
Moreover, to show well-posedness, Stuart \cite{Stuart2010} has proposed a set of sufficient but not necessary assumptions.
Subsequently, Dashti and Stuart \cite{Dashti2017} have reduced these assumptions significantly.
Finally, we mention Ernst et al. \cite{Ernst2015}, who have discussed uniform and H\"older continuity of posterior measures with respect to data, and give sufficient assumptions in this setting. 
We refer to these as \emph{(H\"older, Hellinger)} and \emph{(uniform, Hellinger) well-posedness}, respectively. 

In practical applications, it may be difficult to verify (Lipschitz, Hellinger), (uniform, Hellinger), or (H\"older, Hellinger) well-posedness. The underlying mathematical model can be too complicated to analyse, or even hidden in software.
Indeed, this is the case in large scale applications, e.g., in geotechnical engineering, {meteorology}, and genomics, or in machine learning algorithms.
In any of these cases, the  model is often a black-box -- a function that takes inputs and produces deterministic outputs, but with no known properties.
To the best of our knowledge it is not possible to show (Lipschitz, Hellinger) well-posedness  for the Bayesian inversion of such black-box models.
In turn,  it may not be  necessary to show (Lipschitz, Hellinger) well-posedness for many practical problems; Hadamard's concept contains only continuity, not Lipschitz continuity. In either case, we know that marginal perturbations in the data lead to marginal changes in the posterior measure. Given only continuity, the only difference is that we cannot use information about the data perturbation to quantify the change in the posterior.
This, however, may be tolerable in most practical applications.

Another pressing issue is the measurement of marginal changes in the posterior. Most authors have discussed Lipschitz continuity with respect to the Hellinger distance; exceptions are, e.g.,  the articles of Iglesias et al. \cite{Iglesias2014} and Sprungk \cite{Sprungk2019Lip}.
While the Hellinger distance has useful properties, the actual choice of the metric should depend on the area of application. 
The main contributions of this article are the following:
\begin{enumerate}
\item A new concept of well-posedness of Bayesian inverse problems is proposed. It consists of the existence and the uniqueness of the posterior, as well as of the continuity of the data-to-posterior map in some metric or topological space of probability measures.
\item More specifically, the spaces of probability measures metrised with the Hellinger distance, the total variation distance, and the Wasserstein$(p)$ distance, as well as those associated with the weak topology and the topology induced by the Kullback--Leibler divergence are investigated.
\item Well-posedness of large classes of Bayesian inverse problems in any of these settings is shown. The sufficient assumptions are either non-restrictive or easily verifiable in practical problems (e.g., when having an arbitrary model, finite-dimensional data, and non-degenerate Gaussian noise). The only actually restrictive case remains that of the Kullback--Leibler topology.
\end{enumerate}

This work is organised as follows.
We review the Bayesian approach to inverse problems and the concept of (Lipschitz, Hellinger) well-posedness in \Cref{Sec_IPandIllposed}. 
In \Cref{Sec_WeakWellp}, we advocate our relaxation of Lipschitz continuity and our consideration of metrics other than the Hellinger distance. In the same section, we introduce our notion of well-posedness and show well-posedness with respect to Hellinger, total variation, weak convergence, and the Wasserstein$(p)$ distance, respectively.
In \Cref{Subsec_KLD}, we  extend our concept to stability measurements in the Kullback--Leibler divergence, which is a quasisemi-metric.
We specifically consider the case of finite-dimensional data and non-degenerate Gaussian noise in \Cref{Subsec_Gaussiancase}.
We illustrate our results numerically in  \Cref{Sec_Num_illust} and conclude our work in \Cref{Sec_Conclusions}.
Finally, we note that we review basics on conditional probability in    \Cref{Subsec_Appedix_CondProb} and that we give detailed proofs of all statements formulated in this work in \Cref{Subsec_Appedix_Proofs}.

\section{Inverse problems and the Bayesian approach} \label{Sec_IPandIllposed}
\subsection{Inverse problem}  \label{SubSec_IPs}
Let $y^\dagger$ be observational data in some separable Banach space $(Y, \|\cdot\|_Y)$ -- the \emph{data space}.
The data shall be used to train a mathematical model that is: identify a model parameter $\theta^\dagger$ in a set $X$. The \emph{parameter space} $X$ is a measurable subset of some Radon space $(X', \mathcal{T}')$, i.e. $(X', \mathcal{T}')$ is a separable, completely metrisable topological vector space. $X'$ could, for instance, also be a separable Banach space.
Moreover,  $X$ and $Y$ form measurable spaces with their respective Borel-$\sigma$-algebras $\mathcal{B}X := \mathcal{B}(X,X\cap \mathcal{T}')$ and $\mathcal{B}Y := \mathcal{B}(Y,\|\cdot\|_Y)$.
Let $\mathcal{G}:X \rightarrow Y$ be a measurable function called \emph{forward response operator}. It represents the connection between parameter and data in the mathematical model.
We define the inverse problem by
\begin{equation}\tag{IP}\label{EQ_IP}
\text{Find } \theta^\dagger \in X, \text{ such that } y^\dagger = \mathcal{G}(\theta^\dagger) + \eta^\dagger.
\end{equation}
Here, $\eta^\dagger \in Y$ is observational noise.
We discuss the solvability and stability of inverse problems in terms of well-posedness. 
\begin{definition}[Well-posedness] \label{Def_Wellp_IP}
The problem \cref{EQ_IP} is \emph{well-posed}, if
\begin{enumerate}
\item this problem has a solution, \emph{(existence)}
\item the solution is unique \emph{(uniqueness)}, and
\item the solution depends continuously on the data $y$.  \emph{(stability)}
\end{enumerate}
A problem that is not well-posed is called \emph{ill-posed}.
\end{definition}
We generally consider the observational noise $\eta^\dagger$ to be unknown and model it as a realisation of a random variable $\eta \sim \munoise$.
If the noise takes any value in  $Y$,  the problem \cref{EQ_IP} is ill-posed.
\begin{proposition} \label{prop_ill-posed_IP}
Let $X$ contain at least two elements, and let the support
of $\mu_{\rm noise}$ be $Y$.
Then, the inverse problem \cref{EQ_IP} is {ill-posed}. 
\end{proposition}

Note that the assumptions in \Cref{prop_ill-posed_IP} can often be verified.
If $X$ contains only one element, the inverse problem would be uniquely solvable. However, there is only one possible parameter $\theta \in X$, which makes the inverse problem trivial.
If $Y$ is finite dimensional, the second assumption would for instance be fulfilled, when $\munoise$ is non-degenerate Gaussian.

\subsection{Bayesian inverse problem} \label{subsec_BIP}
The Bayesian approach to \cref{EQ_IP} proceeds as follows:
First, we model the parameter $\theta \sim \muprio$ as a random variable. This random variable reflects the uncertainty in the parameter.
$\muprio$ is the so-called \emph{prior measure}.
Moreover, we assume that $\theta, \eta$ are independent random variables defined on an underlying probability space $(\Omega, \mathcal{A}, \mathbb{P})$.
In this setting, the inverse problem \cref{EQ_IP} is an event:
\begin{equation*}
\{y^\dagger = \mathcal{G}(\theta) + \eta\} \in \mathcal{A},
\end{equation*}
where the data $y^\dagger$ is a realisation of the random variable $\mathcal{G}(\theta^\dagger) + \eta$.
The solution to the Bayesian inverse problem is the posterior measure
\begin{equation} \label{EQ_Post}
\mupost^\dagger := \mathbb{P}(\theta \in \cdot | \mathcal{G}(\theta) + \eta = y^\dagger). 
\end{equation}
For the definition, existence, and uniqueness statement concerning conditional probabilities, we refer to \Cref{Subsec_Appedix_CondProb}.

First, note that we define $y:= \mathcal{G}(\theta) + \eta$ to be a random variable reflecting the distribution of the data, given an uncertain parameter.
We can deduce  the conditional measure of the data $y$ given $\theta = \theta'$: $$\mu_L = \mathbb{P}(y  \in \cdot | \theta = \theta') = \munoise(\cdot -\mathcal{G}(\theta')).$$
To this end, note that the inverse problem setting $y^\dagger = \mathcal{G}(\theta) + \eta$ is only a specific example.
In the following, we consider more general {Bayesian inverse problems}. 
Now, $y$ is a random variable on $(Y, \mathcal{B}Y)$ depending on $\theta$.
The conditional probability of $y$, given that $\theta = \theta'$
is defined by $\mu_L$, which now fully describes the dependence of $\theta$ and $y$. 
The forward response operator $\mathcal{G}$ is implicitly part of $\mu_L$.
\begin{remark}
From a statistical viewpoint, we consider a parametric statistical model $(Y, \mathcal{P})$, where $\mathcal{P} := \{\mu_L(\cdot | \theta') : \theta' \in X \}.$
Hence, the data $y^\dagger \in Y$ is a realisation of $y \sim \mu_L(\cdot|\theta^\dagger)$, for  $\theta^\dagger \in X$. The data $y^\dagger$ is then used to identify this $\theta^\dagger$ among the other elements of $X$. For a thorough discussion of statistical models, we refer to \cite{mccullagh2002}.
\end{remark}
Given $\muprio$ and $\mu_L$, we apply Bayes' Theorem to find the posterior measure $\mupost^\dagger$; now given by
\begin{equation} \label{EQ_Post_gen}
\mupost^\dagger := \mathbb{P}(\theta \in \cdot | y = y^\dagger).
\end{equation}
Bayes' Theorem gives a connection of  $\muprio, \mupost^\dagger$, and $\mu_L$ in terms of their \emph{probability density functions (pdfs)}.
We obtain these pdfs by assuming that there are $\sigma$-finite measure spaces  $(X, \mathcal{B}X, \nu_X)$ and $(Y, \mathcal{B}Y, \nu_Y)$, where $\muprio \ll \nu_X$ and  $\mu_L(\cdot|\theta') \ll \nu_Y$ ($\theta' \in X$, $\muprio$-almost surely (a.s.)).
The Radon--Nikodym Theorem implies that the following pdfs exist:
$$\frac{\mathrm{d}\mu_L}{\mathrm{d} \nu_Y}(y^\dagger) =: L(y^\dagger|\theta'),  \qquad \frac{\mathrm{d}\muprio}{\mathrm{d} \nu_X}(\theta) =: \piprio(\theta).$$
The conditional density $L(\cdot|\theta')$ is called \emph{(data) likelihood}.
The dominating measures $\nu_X$, $\nu_Y$, are often (but not exclusively) given by the counting measure, the Lebesgue measure, or a Gaussian measure. For example, if $X$ is infinite-dimensional and $\muprio$ is Gaussian, we set $\nu_X := \muprio$ and $\piprio \equiv 1$. The posterior measure is then given in terms of a probability density function with respect to the Gaussian prior measure. 
This setting is thorougly discussed in \cite{Dashti2017,Stuart2010}, however, it is also contained in our version of Bayes' Theorem.
Before moving on to that, we discuss a measure-theoretic sublety we encounter with conditional probabilities and their densities.

\begin{remark} \label{Rema_a.s.data}
Conditional probabilities like $\mupost^\dagger= \mathbb{P}(\theta \in \cdot | y = y^\dagger)$ are  uniquely defined  only for $\mathbb{P}(y \in \cdot)$-almost every (a.e.) $y^\dagger \in Y$; see \Cref{Thm_RegCondProb}. 
This implies that if $\mathbb{P}(y \in \cdot)$ has a continuous distribution, point evaluations in $Y$ of the function $\mathbb{P}(\theta \in A | y = \cdot)$ may not be well-defined, for $A \in \mathcal{B}X$.
In this case, one would not be able to compute the posterior measure for any single-point data set $y^\dagger \in Y$.
Also, the statements \cref{EQ_Post}, \cref{EQ_Post_gen}, as well as the definition of the likelihood, should be understood only for $\mathbb{P}(y \in \cdot)$-a.e. $y^\dagger \in Y$. 
\end{remark}
Our version of Bayes' Theorem is mainly built on  \cite[Theorem 3.4]{Dashti2017}. However, in the proof we neither need to assume that the model evidence is positive and finite, nor do we need to assume continuity in the data or the parameter of the likelihood. 

\begin{theorem}[Bayes]\label{thm:bayes} Let $y^\dagger \in Y$ be $\mathbb{P}(y\in \cdot)$-almost surely defined. Moreover, let $L(y^\dagger|\cdot)$ be in $\mathbf{L}^1(X,\muprio)$ and  strictly positive. 
Then, $$Z(y^\dagger) := \int L(y^\dagger|\theta) \mathrm{d}\muprio(\theta) \in (0, \infty).$$ Moreover, the posterior measure $\mupost^\dagger \ll \nu_X$ exists, it is unique, and it has the $\nu_X$-density
\begin{equation} \label{EQ_Bayes}
 \pipost^\dagger(\theta') = \frac{L(y^\dagger|\theta')\piprio(\theta')}{Z(y^\dagger)} \qquad (\theta' \in X, \nu_X\text{-a.s.}).
\end{equation}
\end{theorem}


The quantity in the denominator of Bayes' formula $Z(y^\dagger) := \int L(y^\dagger|\theta) \mathrm{d}\muprio(\theta)$ is the $\nu_Y$-density of $\mathbb{P}(y \in \cdot)$ and is called \emph{(model) evidence}.
We comment on the assumptions of \Cref{thm:bayes} in \Cref{SubseHellTVWeak}.
In \Cref{Rema_a.s.data}, we mention that the posterior measure is only $\mathbb{P}(y \in \cdot)$-a.s. uniquely defined. Hence, the map $y^\dagger \mapsto \mupost^\dagger$ is not well-defined. We resolve this issue by fixing the definition of the likelihood $L(y^\dagger|\theta')$ for every $y^\dagger \in Y$ and $\muprio$-a.e. $\theta' \in X$.
According to \Cref{thm:bayes}, we then obtain indeed a unique posterior measure for any data set $y^\dagger \in Y$.
We define the \emph{Bayesian inverse problem} with prior $\muprio$ and likelihood $L$ by
\begin{equation}\tag{BIP}\label{EQ_BIP}
\text{Find } \mupost^\dagger\in \Prob(X, \muprio) \text{ with $\nu_X$-density } \pipost(\theta|y^\dagger)= \frac{L(y^\dagger|\theta)\piprio(\theta)}{Z(y^\dagger)}.
\end{equation}
Here,  $\Prob(X,\muprio)$ denotes the set of probability measures on $(X, \mathcal{B}X)$ which are absolutely continuous with respect to the prior $\muprio$.
Similarly, we define the set of all probability measures  on $(X, \mathcal{B}X)$ by $\Prob(X)$. If $X$ forms a normed space with some norm $\|\cdot\|_X$, we define the set of probability measures with finite $p$-th moment by
$$\Prob_p(X) := \left\lbrace \mu \in \Prob(X): \int \|\theta\|_X^p \mu(\mathrm{d}\theta)  < \infty \right\rbrace  \qquad (p \in [1, \infty)).$$

\subsection{Degenerate Bayesian inverse problems} \label{SubSec_NoBayes}
There are Bayesian inverse problems, for which Bayes' Theorem (\Cref{thm:bayes}) is not satisfied.  Consider $\munoise := \delta(\cdot - 0)$ as a noise distribution; i.e. the noise is almost surely $0$. We refer to Bayesian inverse problems with such a noise distribution as \emph{degenerate}, since the noise distribution is {degenerate}.
Here, we represent the likelihood by
$$
L(y^\dagger|\theta') := \begin{cases} 1, &\text{ if } y^\dagger = \mathcal{G}(\theta'),\\
0,  &\text{ otherwise.}
\end{cases}
$$
Due to different dimensionality, it is now likely that the prior $\muprio$ is chosen such that it gives probability $0$ to the solution manifold
$
S = \{ \theta' \in X : y^\dagger = \mathcal{G}(\theta')\},
$
i.e. $\muprio(S) = 0$. Then, we have $$Z(y^\dagger) =\int_X L(y^\dagger|\theta) \muprio(\mathrm{d}\theta) = \int_S 1 \muprio(\mathrm{d}\theta) = \muprio(S) = 0$$ and do not obtain a valid posterior measure for $y^\dagger$ from \Cref{thm:bayes}.
Alternatively, we can employ the Disintegration Theorem, see Cockayne et al. \cite{Cockayne2017} and the definition of conditional probabilities after \Cref{Thm_RegCondProb}.
In the following proposition, we give a simple example for such a \cref{EQ_BIP}.
\begin{proposition} \label{prop_homeomo}
Let $\mathcal{G}: X \rightarrow Y$ be a homeomorphism, i.e. it is continuous, bijective, and $\mathcal{G}^{-1}: Y \rightarrow X$ is continuous as well.  Moreover, let $\muprio \in \Prob(X)$ be some prior measure and $\munoise = \delta( \cdot - 0)$.
Then, 
$\mupost^\dagger = \delta(\mathcal{G}(\cdot) - y^\dagger) =  \delta(\cdot - \mathcal{G}^{-1}(y^\dagger)), $ for $\muprio(\mathcal{G} \in \cdot)$-a.e. $y^\dagger \in Y$.
\end{proposition}
Note that we cannot easily solve the problem discussed in \ref{Rema_a.s.data} for this Bayesian inverse problem. Hence, point evaluations $y^\dagger \mapsto \mupost^\dagger$ may indeed be not well-defined in this setting. 
Therefore, when discussing this Bayesian inverse problem, we will fix one representative in the class of measures that are almost surely equal to the posterior.
\subsection{Lipschitz well-posedness}
We now move on to the definition of  Stuart's \cite{Stuart2010} concept of well-posedness of Bayesian inverse problems.
Similarly to the well-posedness definition of the classical problem \cref{EQ_IP}, we consider an existence, a uniqueness and a stability condition; see \Cref{Def_Wellp_IP}.
Stability is quantified in terms of the 
\emph{Hellinger distance}
$$\dHel(\mu, \mu') = \sqrt{\frac{1}{2}\int  \left(\sqrt{\frac{\mathrm{d} \mu'}{\mathrm{d}\muprio}} - \sqrt{\frac{\mathrm{d} \mu}{\mathrm{d}\muprio}} \right)^2 \mathrm{d}\muprio}$$
between two measures $\mu, \mu' \in \Prob(X,\muprio)$.
The Hellinger distance is based on the work \cite{Hellinger1909}. With this, we can now formalise the concept of (Lipschitz, Hellinger) well-posedness for Bayesian inverse problems.
\begin{definition}[(Lipschitz, Hellinger) well-posedness] \label{Def_well-posed_BIP}
The problem \cref{EQ_BIP} is \emph{(Lipschitz, Hellinger) well-posed}, if
\begin{enumerate}
\item $\mupost^\dagger \in \Prob(X, \muprio)$ exists, \emph{(existence)}
\item $\mupost^\dagger$ is unique in $\Prob(X, \muprio)$ \emph{(uniqueness)}, and 
\item $(Y, \|\cdot\|_Y) \ni y^\dagger \mapsto \mupost^\dagger \in (\Prob(X, \muprio),\dHel)$ is locally Lipschitz continuous. \emph{(stability)}
\end{enumerate}
\end{definition}

\section{Redefining Well-posedness} \label{Sec_WeakWellp}

In this work, we try to identify general settings in which we can show some kind of well-posedness of \cref{EQ_BIP}, using no or very limited assumptions on the underlying mathematical model or the forward response operator.
In particular, we aim to find verifiable assumptions on the likelihood  $L(y^\dagger|\theta')$ (or rather the noise model) that are independent of the underlying forward response operator $$\mathcal{G} \in \mathbf{M}:=\{f:X\rightarrow Y \text{ measurable}\}.$$

Neglecting \Cref{prop_homeomo} for a moment, existence and uniqueness are often implied by \Cref{thm:bayes}. 
However, the local Lipschitz continuity condition, reflecting stability, is rather strong.
In \Cref{Subsec_Motiv}, we give examples in which local Lipschitz continuity does not hold in the posterior measure or is hard to verify by using results in the literature.
In any of these cases, we show that the posterior measures are continuous in the data.
Given that the classical formulation of well-posedness, i.e. \Cref{Def_Wellp_IP},  does not require local Lipschitz continuity and that local Lipschitz continuity may be too strong for general statements, we use these examples to advocate a relaxation of the local Lipschitz continuity condition.

Moreover, it is not possible to use the Hellinger distance to quantify the distance between two posteriors  in some situtations. In other situations, the Hellinger distance may be inappropriate from a contentual viewpoint. In \Cref{Subsec_Recon_Hell}, we will investigate these issues as a motivation to consider metrics other than the Hellinger distance.

In \Cref{Subsec_Def_Main_res}, we will introduce the concept of $(\PP, d)$-wellposedness of Bayesian inverse problems. Finally, we will show the main results of this work: we give conditions under which we can show well-posedness in a variety of metrics in \Cref{SubseHellTVWeak,SubsecWass}.
\subsection{Relaxing the Lipschitz condition} \label{Subsec_Motiv}
\begin{figure}
\centering
\includegraphics[scale=0.69]{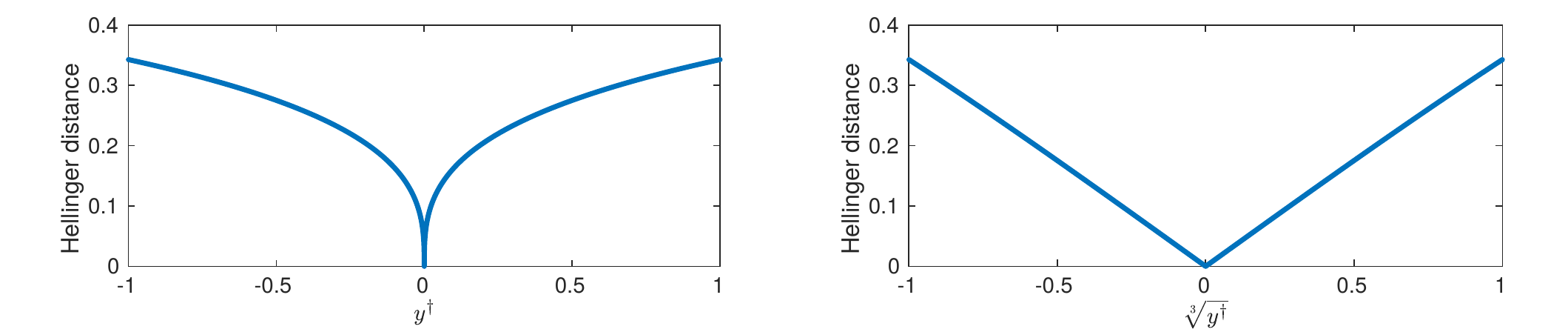}
\caption{Hellinger distances between posterior measures in \Cref{Ex_cubic_text}. The posterior measures are based on two data sets: $y^\dagger$ that varies in (-1,1) and $y^\ddagger := 0$. 
In the left figure, we show the relationship between data and Hellinger distance. In the right figure, we replace the data by $y^\dagger := \sqrt[3]{y^\dagger}$, $y^\ddagger := \sqrt[3]{y^\ddagger}$. In both plots, we observe a continuous relationship between Hellinger distance and data, which is also Lipschitz continuous in the right figure, but not in the left figure.}
\label{Fig_Hell3rdroot}
\end{figure}

Ill-posedness in the (Lipschitz, Hellinger) sense can for instance occur when data has been transformed by a non-Lipschitz continuous function. 
As an example, we consider a Bayesian inverse problem that is  linear and Gaussian, however, the data is transformed by the cube root function.
\begin{example} \label{Ex_cubic_text}
Let $X:= Y:= \mathbb{R}$. 
We consider the Bayesian approach to the inverse problem
$$
y^\dagger = (\theta + \eta)^{3},
$$
where $\theta$ is the unknown parameter and $\eta$ is observational noise; both are independent. The probability measure of parameter and noise are given by $\muprio := \munoise := \mathrm{N}(0,1^2)$. 
The likelihood of the  \cref{EQ_BIP} is
$$
L(\theta|y^\dagger) = \frac{1}{\sqrt{2\pi}} \exp\left(-\frac{1}{2}\|\theta - \sqrt[3]{y^\dagger}\|^2\right).
$$
Since prior and noise are Gaussian, and the forward model is linear (the identity operator), we can compute the posterior measure analytically, see \cite[\S 3]{Agapiou14}. We obtain $\mupost^\dagger := \mathrm{N}\left(\sqrt[3]{y^\dagger}/2,({1}/{\sqrt{2}})^2 \right)$.
Moreover, one can show that 
\begin{align} \label{EQ_Hellinger_cubic}
\dHel(\mupost^\dagger, \mupost^\ddagger) &= \sqrt{1- \exp\left(-\frac18 \left(\sqrt[3]{y^\dagger}-\sqrt[3]{y^\ddagger}\right)^2\right)},
\end{align}
where $\mupost^\ddagger$ is the posterior measure based on a second data set $y^\ddagger \in Y$.
One can show analytically that this Hellinger distance in \cref{EQ_Hellinger_cubic} is not locally Lipschitz as $|y^\dagger - y^\ddagger| \rightarrow 0$. It is however continuous. 
We plot the Hellinger distance in \Cref{Fig_Hell3rdroot} on the left-hand side, where we set $y^\ddagger := 0$ and vary only $y^\dagger \in (-1,1)$.
We observe indeed  that the Hellinger distance is continuous, but not Lipschitz continuous.
In the plot on the right-hand side, we show the Hellinger distance, when considering $\sqrt[3]{y^\dagger}$ as the data set, rather than $y^\dagger$.
In this case, the Hellinger distance is locally Lipschitz in the data.
\end{example}
The Bayesian inverse problem in \Cref{Ex_cubic_text} is  ill-posed in the sense of \Cref{Def_well-posed_BIP}, since the posterior is only continuous but not Lipschitz in the data.
However, we can heal this ill-posedness by transforming $y^\dagger \mapsto \sqrt[3]{y^\dagger}$.
Hence, the (Lipschitz, Hellinger) well-posedness property reduces to a continuous data transformation problem.

Other examples may be (Lipschitz, Hellinger) well-posed, but this may be difficult to verify in practice, or for general forward response operators.
Dashti and Stuart give  \cite[Assumptions 1]{Dashti2017} that are sufficient, but not necessary, to prove well-posedness.
One of the assumptions is local Lipschitz continuity in the log-likelihood $\log L$ with respect to the data. 
Here, the Lipschitz constant is supposed to be a positive function that is monotonically non-decreasing in $\|\theta\|_X$.
This assumption is not satisfied in the following example.
\begin{example} \label{Ex_inverse_text}
Let $X:= (0,1)$ and $Y:= \mathbb{R}$. 
We consider the Bayesian approach to the inverse problem
$$
y^\dagger = \theta^{-1} + \eta
$$
where $\theta$ is the unknown parameter and $\eta$ is observational noise.
Neglecting linear prefactors, this inverse problem can be thought of as the recovery of a wavelength $\theta$ from a noisy frequecy measurement $y^\dagger$. 

The prior measure of $\theta$ is given by $\muprio = \mathrm{Unif}(0,1)$. The noise is distributed according to $\munoise = \mathrm{N}(0,1^2)$. 
Moreover, note that parameter and noise are independent random variables.
The likelihood of the  \cref{EQ_BIP} is
$$
L(y^\dagger|\theta) = \frac{1}{\sqrt{2\pi}} \exp\left(-\frac{1}{2}\|\theta^{-1} - {y^\dagger}\|^2\right).
$$
For fixed $\theta \in X$, the logarithm of the likelihood in this setting is Lipschitz continuous in the data. 
However, as $\theta \downarrow 0$, the Lipschitz constant explodes.
Hence, the likelihood does not fulfil \cite[Assumptions 1]{Dashti2017}.
\end{example}
Hence, we cannot use the theory of Dashti and Stuart \cite[\S 4]{Dashti2017} to show (Lipschitz, Hellinger) well-posedness of the Bayesian inverse problem in \Cref{Ex_inverse_text}.
We expect a similar problem for forward response operators that are not locally bounded.
In \Cref{Cor_Gaussian_case}, we  revisit \Cref{Ex_inverse_text} and show that the posterior measure is continuous with respect to the data. 

Up to now we presented academic examples.
A practically more relevant problem is the Bayesian elliptic inverse problem.
It is the prototype example in the context of partial differential equations and has been investigated by various authors, e.g., \cite{Dashti2011,Dashti2017,LEU19,Richter1981, Stuart2010}. 
\begin{example}[Elliptic inverse problem] \label{Ex_Elliptic}
Let the parameter space be a space of continuous functions $X := \mathbf{C}^0(D)$ on a bounded open set $D \subseteq \mathbb{R}^d$, $d=1,2,3$. The data space $Y:= \mathbb{R}^k$ is finite-dimensional.
The underlying model is an elliptic partial differential equation:
\begin{align*}
-\nabla \cdot \left(\mathrm{e}^{\theta(x)}\,\nabla p(x)\right) &= f(x) &(x \in D) \\
p(x) &= 0 &(x \in \partial D),
\end{align*}
or rather its weak formulation.
In a typical application, the solution $p$ represents the pressure head in a groundwater reservoir, while the diffusion coefficient $\exp(\theta(x))$ represents the reservoir's hydraulic conductivity.
Noisy measurements of the pressure head at locations $x_1,\ldots,x_k \in D$ shall now be used to infer the log-conductivity $\theta$. 
Hence, the forward response operator is the map
$$
\mathcal{G}:X\rightarrow Y, \qquad \theta \mapsto (p(x_1),\ldots,p(x_k)).
$$
\end{example}
In practical applications, allowing only continuous functions as diffusion coefficients may be too restrictive. 
Iglesias et al. \cite{Iglesias2014} consider more realistic geometric priors measures. In \cite[Theorem 3.5]{Iglesias2014}, the authors show local Lipschitz continuity for some of those prior measures, but only H\"older continuity with coefficient $\gamma = 0.5$ for others. 
This is another example where (Lipschitz, Hellinger) well-posedness in the sense of \Cref{Def_well-posed_BIP} has not been shown, but continuity in the posterior measure is satisfied.

In the following, in \Cref{Subsec_Def_Main_res}, we weaken the \emph{Lipschitz well-posedness} by replacing Lipschitz continuity with continuity as a stability condition.
Looking back at \Cref{Ex_cubic_text,Ex_inverse_text,Ex_Elliptic}, we consider this weakening {tolerable} for practical problems.

\subsection{Reconsidering the Hellinger distance} \label{Subsec_Recon_Hell}
The Hellinger distance is a popular choice to analyse the continuous dependence or, e.g., the approximation of measures.
However, there are cases, in which it cannot be used:

We consider the Bayesian inverse problem discussed in \Cref{prop_homeomo}. We set $\mupost^\dagger := \delta(\cdot - \mathcal{G}^{-1}(y^\dagger))$ as a posterior measure with $\mathcal{G}^{-1}: Y \rightarrow X$ continuous.
We set $X := Y := \mathbb{R}$ and $\muprio := \mathrm{N}(0,1^2)$. Then, $\mupost^\dagger \not\ll \muprio$, for $y^\dagger \in Y$.
The Hellinger distance between $\mupost^\dagger$ and $\mupost^\ddagger$ is not well-defined, for any other data set $y^\ddagger \neq y^\dagger$.
Instead, we consider the closely related total variation (tv) distance and obtain
$$ \dtv(\mupost^\dagger, \mupost^\ddagger):= \sup_{B \in \mathcal{B}X}\Big|\mupost^\dagger(B)- \mupost^\ddagger(B)\Big| =1.$$
Hence, $\mupost^\ddagger \not\rightarrow \mupost^\dagger$ in total variation as $y^\ddagger \rightarrow y^\dagger$.
Thus, the Bayesian inverse problem is not stable in the total variation distance, i.e. ill-posed in this sense.

However, we have $\mupost^\ddagger \rightarrow \mupost^\dagger$ weakly as $y^\ddagger \rightarrow y^\dagger$. 
Hence, we observe continuity in the weak topology  on the space $\Prob(X)$ of probability measures on $(X, \mathcal{B}X)$.
Equivalently, we can say that we observe continuity in the \emph{Prokhorov metric} on $\Prob(X)$:
$$ \dProk(\mu, \mu') := \inf\left\lbrace \varepsilon > 0: \mu(B) \leq \mu'(B^\varepsilon) + \varepsilon, B \in \mathcal{B}X \right\rbrace,$$
where $B^\varepsilon := \{b \in X : b' \in B, \|b-b'\|_X < \varepsilon\}$ is the open generalised $\varepsilon$-ball around $B \in \mathcal{B}X$; see \cite{Prokhorov56} for details.

To summarise, there are cases in which the Hellinger distance is infeasible to show well-posedness. Moreover, different metrics may lead to different well-posedness results. 
Hence, we should introduce a concept that allows for different metrics on the space of probability measures.

\subsection{Definition} \label{Subsec_Def_Main_res}
As motivated in \Cref{Subsec_Motiv,Subsec_Recon_Hell}, we next generalise the notion of well-posedness of Bayesian inverse problems. In \Cref{Def_well-posed_BIP}, we considered Lipschitz continuity in the Hellinger distance as a stability criterion. 
Now, we consider simple continuity with respect to various metric spaces.
\begin{definition}[$(\PP,d)$-Well-posedness] \label{Def_well-posed_BIP_NEW}  Let $(\PP,d)$ be a metric space of probability measures on $(X, \mathcal{B}X)$; i.e. $\PP \subseteq \Prob(X)$. The problem \cref{EQ_BIP}
 is $(\PP,d)$-\emph{well-posed}, if
\begin{enumerate}
\item $\mupost^\dagger \in \PP$ exists, \emph{(existence)}
\item $\mupost^\dagger$ is unique in $\PP$ \emph{(uniqueness)}, and
\item $(Y, \|\cdot\|_Y) \ni y^\dagger \mapsto \mupost^\dagger \in (\PP,d)$ is a  continuous function. \emph{(stability)}
\end{enumerate}
\end{definition}
For particular $(\PP, d)$, we introduce special denominations.
Indeed, 
we denote $(\PP, d)$-well-posedness by
\begin{itemize}
\item[(i)] \emph{weak} well-posedness, if we consider the \emph{Prokhorov metric}, i.e. we set $(\PP, d)=(\Prob(X), \dProk)$,
\item[(ii)]
 \emph{total variation} well-posedness, if we consider the \emph{total variation distance}, i.e.  we set $(\PP, d):=(\Prob(X), \dtv)$, 
 \item[(iii)]  \emph{Hellinger} well-posedness, if we consider the  \emph{Hellinger distance}, i.e. we set $(\PP, d):=(\Prob(X,\muprio), \dHel)$,  and
\item[(iv)]  \emph{Wasserstein($p$)} well-posedness, if $X$ is a normed space and if we consider the \emph{Wasserstein($p$) distance}, i.e. we set $(\PP, d):=(\Prob_p(X), \dWas)$, for some $p \in [1,\infty)$. 
\end{itemize}
Which concept of  well-posedness should we consider in practice?
Weak well-posedness implies continuity of posterior expectations of bounded, continuous quantities of interest.
If this is the task of interest, weak well-posedness should be sufficient. Hellinger and tv distance imply convergence of the posterior expectation of any bounded quantity of interest. Hence, if discontinuous functions shall be integrated, or probabilities computed, those distances should be chosen. 
Wasserstein$(p)$ distances have gained popularity in the convergence and stability theory of Markov chain Monte Carlo (MCMC) algorithms, see, e.g., \cite{Gibbs2004,rudolf2018}.
Hence,  Wasserstein$(p)$ well-posedness may be the right tool when discussing the well-posedness of solving a Bayesian inverse problem via MCMC.

%
\subsection{Hellinger, total variation, and weak well-posedness} \label{SubseHellTVWeak}
We now give assumptions, under which a Bayesian inverse problem can be shown to be Hellinger well-posed, total variation well-posed, and weakly well-posed. 
\begin{assumptions} \label{Assumptions_main} Consider a \cref{EQ_BIP}.
Let the following assumptions hold for $\muprio$-almost every $\theta' \in X$ and every  $y^\dagger \in Y$.
\begin{itemize}
\item[(A1)] $L(\cdot |\theta')$ is a strictly positive probability density function,
\item[(A2)] $L(y^\dagger |\cdot) \in \mathbf{L}^1(X,\muprio)$,
\item[(A3)] $g \in \mathbf{L}^1(X,\muprio)$ exists such that $L(y^\ddagger |\cdot) \leq g$ for all $y^\ddagger \in Y$, and
\item[(A4)]  $L(\cdot|\theta')$ is continuous.
\end{itemize}
\end{assumptions}
(A1) means that any data set $y^\dagger \in Y$ has a positive likelihood under any parameter $\theta' \in X$. We conservatively assume that no combination of parameter and data values is impossible, but some may be unlikely. 
This can usually be satisfied by continuously transforming the forward response operator and/or by choosing a noise distribution that is concentrated on all of $Y$.
Note that the assumption that $L(y^\dagger|\theta')$ is a probability density function can be relaxed to $c \cdot L(y^\dagger|\theta')$ is a probability density function, where $c > 0$ does depend neither on $y^\dagger$, nor $\theta'$.
(A2)-(A3) imply that the likelihood is integrable with respect to the prior and that it is bounded from above uniformly in the data by an integrable function.
These assumptions are for instance satisfied, when the likelihood is bounded from above by a constant.
Noise models with bounded probability density function on $Y$ should generally imply a bounded likelihood.
(A4) requires the continuity of the likelihood with respect to the data.
Continuity in the data is for instance given, when considering noise models with continuous probability density functions, and a continuous  connection of noise and model.
We give examples in \Cref{Sec_Num_illust} showing that we cannot neglect the continuity in the data.
Here, we show Hellinger, total variation, and weak well-posedness under (A1)-(A4).

\begin{theorem} \label{Thm_main}
Let (A1)-(A4) hold for a \cref{EQ_BIP}. Then, \cref{EQ_BIP} is weakly, Hellinger, and total variation (whtv) well-posed.
\end{theorem}

For the proof of this theorem, we proceed as follows: first, we show that (A1)-(A4) imply Hellinger well-posedness. Then we show that total variation well-posedness and weak well-posedness are indeed implied by Hellinger well-posedness by some topological argument. 

\begin{lemma} \label{Lemma_main}
Let (A1)-(A4) hold for a \cref{EQ_BIP}. Then, \cref{EQ_BIP} is Hellinger well-posed.
\end{lemma}
We can bound Prokhorov and total variation distance with the Hellinger distance; see  \cite{Gibbs2002} for the appropriate results. 
In such a case, the continuity of a function in the bounding metric immediately implies continuity also in the bounded metric. 
\begin{lemma} \label{Lem_topologies}
Let $A, B$ be two sets and let $(A, d_A)$, $(B, d_1)$ and $(B, d_2)$ be metric spaces.
Let $f: (A, d_A) \rightarrow (B, d_2)$ be a continuous function. Moreover, let $t: [0, \infty) \rightarrow [0, \infty)$ be  continuous in $0$, with $t(0)  = 0$.
Finally, let be
$$ d_1(b,b') \leq t(d_2(b,b')) \qquad (b, b' \in B).$$
Then, $f: (A, d_A) \rightarrow (B, d_1)$ is continuous as well.
\end{lemma}
In the setting of \Cref{Lem_topologies}, we call $d_1$ \emph{coarser} than $d_2$, respectively  $d_2$ \emph{finer} than $d_1$.
The lemma implies that if we are going from a finer to a coarser metric, continuous functions keep on being continuous.
By the bounds given in \cite{Gibbs2002}, Prokhorov and total variation distance are coarser than the Hellinger distance. 
If the function $y^\dagger \mapsto \mupost^\dagger$ is continuous in the Hellinger distance, it is also continuous in the weak topology and the total variation distance. 
We summarise this result in the following proposition.
\begin{proposition} \label{Propo_wellposedd1d2}
Let $d_1, d_2$ be metrics on $\PP$ and let $d_1$ be coarser than $d_2$.
Then, a Bayesian inverse problem that is  $(\PP, d_2)$-wellposed,  is also   $(\PP, d_1)$-wellposed.
\end{proposition}
Therefore, Hellinger well-posedness (in \Cref{Lemma_main}) implies total variation and weak well-posedness (in \Cref{Thm_main}).
\subsection{Wasserstein$(p)$ well-posedness} \label{SubsecWass}
Let $p \in [1,\infty)$ and let $X$ form a normed space with norm $\|\cdot \|_X$. The \emph{Wasserstein$(p)$ distance} between $\mu, \mu' \in \Prob_p(X)$ can be motivated by the theory of optimal transport. It is given as the cost of the optimal transport from $\mu$ to $\mu'$. The cost of transport from $\theta \in X$ to $\theta' \in X$ is given by $\|\theta - \theta'\|_X$.
More precisely, the \emph{Wasserstein$(p)$} distance (i.e. the Wasserstein distance of order $p$) is defined by
$$\dWas(\mu, \mu') := \left(\inf_{\Lambda \in \mathrm{C}(\mu, \mu')}\int_{X \times X}\|\theta - \theta'\|_X^p\mathrm{d}\Lambda(\theta,\theta')\right)^{1/p},$$
where
$\mathrm{C}(\mu, \mu') := \lbrace \Lambda' \in \Prob(X^2) : \mu(B) = \Lambda'(B \times X), \mu'(B) = \Lambda'(X \times B), B \in \mathcal{B}X\rbrace$ is the \emph{set of couplings} of $\mu, \mu' \in \Prob_p(X)$.
We can link convergence in the Wasserstein$(p)$ distance to weak convergence.
Let $(\mu_n)_{n \in \mathbb{N}} \in \Prob_p(X)^{\mathbb{N}}$ be a sequence and $\mu \in \Prob_p(X)$ be some other probability measure. Then, according to \cite[Theorem 6.9]{villani2009optimal}, we have
\begin{align}
&\lim_{n\rightarrow \infty}\dWas(\mu_n, \mu) = 0 \label{Eq_WasProk} \\ \Leftrightarrow &\left( \lim_{n\rightarrow \infty}\dProk(\mu_n, \mu) = 0 \text{ and } \lim_{n\rightarrow \infty}\int\|\theta\|_X^p \mu_n(\mathrm{d}\theta) =  \int\|\theta\|_X^p \mu(\mathrm{d}\theta) \right). \nonumber
\end{align}
Hence, to show Wasserstein$(p)$ well-posedness, we need to show weak well-posedness and stability of the $p$-th posterior moment with respect to changes in the data. 
Assumptions (A1)-(A4) are not sufficient to show the latter.
As in \Cref{SubseHellTVWeak}, we now give the additional assumption (A5) that we need to show Wasserstein$(p)$ well-posedness. 
Then, we discuss situations in which this assumption is satisfied. We finish this section by showing Wasserstein well-posedness under (A1)-(A5).
\begin{assumptions} \label{Assumptions_Wasser} Consider a \cref{EQ_BIP}.
Let the following assumption hold.
\begin{itemize}
\item[(A5)] $g' \in \mathbf{L}^1(X,\muprio)$ exists such that $\|\theta'\|_X^p \cdot L(y^\dagger |\theta') \leq g'(\theta')$ for $\muprio$-a.e. $\theta' \in X$ and all $y^\dagger \in Y$.
\end{itemize}
\end{assumptions}

Assumption (A5)  eventually requires a uniform bound on the $p$-th moment of the posterior measure. This is in general not as easily satisfied as (A1)-(A4). However, there is a particular case, in which we can show that (A1)-(A5) are satisfied rather easily: 
If the likelihood is bounded uniformly by a constant and if the prior has a finite $p$-th moment.
\begin{proposition}\label{PropA1A5}
 We consider a \cref{EQ_BIP} and some $p \in [1,\infty)$. Let (A1) and (A4) hold. Moreover, let
 some $c \in (0, \infty)$ exist, such that 
 $$L(y^\dagger|\theta') \leq c \qquad (y^\dagger \in Y; \theta' \in X, \muprio\text{-a.s.}), $$  
 and let $\muprio \in \Prob_p(X)$. 
 Then (A1)-(A5) are satisfied.
\end{proposition}
We have already mentioned that a uniformly bounded likelihood does not appear to be a very restrictive property. 
Boundedness of the $p$-th moment of the prior is rather restrictive though. 
In practical problems, prior measures very often come from known well-known families of probability measures, such as Gaussian, Cauchy, or exponential.  For such families we typically know whether certain moments are finite.
In this case, it is easy to see with \Cref{PropA1A5}, whether the \cref{EQ_BIP} satisfies assumption (A5). 
Hence, (A5) is restrictive, but easily verifiable.
Next, we state our result on Wasserstein$(p)$ well-posedness.

\begin{theorem} \label{thm:WassWell}
Let $p \in [1,\infty)$ and let (A1)-(A5)  hold for a \cref{EQ_BIP}. Then, \cref{EQ_BIP} is Wasserstein$(p)$ well-posed.
\end{theorem}

Finally, we note that weak and Wasserstein$(p)$ stability can also hold in degenerate Bayesian inverse problems; see \Cref{SubSec_NoBayes}. 
Given the argumentation in \Cref{Subsec_Recon_Hell}, we see that the Bayesian inverse problem discussed in \Cref{prop_homeomo} is stable in the weak topology but neither in the Hellinger nor in the total variation sense. 
Indeed, the Bayesian inverse problem is also stable in the Wasserstein$(p)$ distance for any $p \in [1,\infty)$, but neither satisfies (A1), nor (A4).
\begin{corollary} \label{Coro_homeo_wellp}
We consider the Bayesian inverse problem given in \Cref{prop_homeomo}, i.e. we assume that the posterior measure is given by $$\mupost^\dagger = \delta(\cdot - \mathcal{G}^{-1}(y^\dagger)) \qquad (y^\dagger \in Y),$$
and $\mathcal{G}^{-1}: Y \rightarrow X$ is continuous. 
Then, this posterior measure stable in the weak topology. If $X$ is additionally a normed space, the posterior is also stable in Wasserstein$(p)$, for any $p \in [1,\infty)$. 
\end{corollary}

\section{Well-posedness in quasisemi-metrics} \label{Subsec_KLD}
The distances  we have considered in \Cref{Sec_WeakWellp} ($\dHel, \dtv, \dProk, \dWas$) are all well-defined metrics. In statistics and especially in information theory, various distance measures are used that are not actually metrics. For instance, they are asymmetric (quasi-metrics), they do not satisfy the triangle inequality (semi-metrics), or they do not satisfy either (quasisemi-metrics).
Due to their popularity, it is natural to consider stability also in such generalised distance measures. 

The \emph{Kullback--Leibler divergence} (KLD), \emph{relative entropy}, or \emph{directed divergence} is a popular quasisemi-metric used in information theory and machine learning.
In the following, we consider the KLD exemplarily as a quasisemi-metric, in which we discuss well-posedness.
The KLD is used to describe the \emph{information gain} when going from $\mu \in \Prob(X)$ to another measure $\mu' \in \Prob(X,\mu)$. If defined, it is given by
$$\DKL(\mu'\|\mu)  := \int_X \log\left(\frac{\mathrm{d}\mu'}{\mathrm{d}\mu} \right)\mathrm{d}\mu'.$$
The KLD induces a topology; see \cite{Belavkin2015}. Hence, we can indeed describe continuity in the KLD and, thus, consider the \emph{Kullback--Leibler well-posedness} of Bayesian inverse problems.
This concept bridges information theory and Bayesian inverse problems; and allows  statements about the loss of information in the posterior measure, when the data is perturbed. In particular, we define this \emph{loss of information} by the information gain when going from the posterior $\mupost^\ddagger$ with perturbed data $y^\ddagger$ to the posterior $\mupost^\dagger$ with unperturbed data $y^\dagger$. Hence, the loss of information is equal to $\DKL(\mupost^\dagger\|\mupost^\ddagger).$
A Bayesian inverse problem is Kullback--Leibler well-posed, if the posterior measure exists, if it is unique, and if the information loss is continuous with respect to the data.
\begin{definition}[Kullback--Leibler well-posed] \label{Def_KLWell}
The problem \cref{EQ_BIP} is \emph{Kullback--Leibler well-posed}, if
\begin{enumerate}
\item $\mupost^\dagger \in \Prob(X, \muprio)$ exists \emph{(existence)},
\item $\mupost^\dagger$ is unique in $\Prob(X, \muprio)$ \emph{(uniqueness)}, and
\item for all $y^\dagger \in Y$ and $\varepsilon > 0$, there is $\delta(\varepsilon) > 0$, such that  
$$
\DKL(\mupost^\dagger\|\mupost^\ddagger) \leq \varepsilon \quad (y^\ddagger \in Y: \|y^\dagger - y^\ddagger\|_Y \leq \delta(\varepsilon))
\qquad \text{\emph{(stability)}}.$$ 
\end{enumerate}
\end{definition}
In the setting of  \cref{thm:bayes}, 
(A1)-(A4) are not sufficient to show Kullback--Leibler well-posedness; indeed, the Kullback--Leibler divergence may be not even well-defined.
We require the following additional assumption on the log-likelihood. 
\begin{assumptions} \label{Assumptions_KLD} Consider a \cref{EQ_BIP}.
Let the following assumption hold for $\muprio$-almost every $\theta' \in X$ and every  $y^\dagger \in Y$:
\begin{itemize}
\item[(A6)] there is a $\delta > 0$ and a function $h(\cdot,y^\dagger) \in \mathbf{L}^1(X,\mupost^\dagger)$ such that $$|\log L(y^\ddagger|\cdot)| \leq h(\cdot, y^\dagger) \qquad \qquad (y^\ddagger \in Y :  \|y^\dagger - y^\ddagger\|_Y \leq \delta).$$
\end{itemize}
\end{assumptions}
Assumption (A6) is much more restrictive than  (A1)-(A4). 
Indeed, we now require some integrability condition on the forward response operator.
The condition may be hard to verify,  when the posterior measure has heavy tails or when the model is unbounded. Also, when we are not able to analyse the underlying model.

\begin{theorem} \label{thm:KLD}
Let (A1)-(A4), (A6) hold for a \cref{EQ_BIP}. Then, \cref{EQ_BIP} is Kullback--Leibler well-posed.
\end{theorem}
\begin{remark}\label{Remark_local_versions_A}
We note that we have allowed the bound in (A6) to depend on $y^\dagger \in Y$ and to hold only locally on sets of the form $\{y^\ddagger : \|y^\dagger - y^\ddagger\|_Y \leq \delta\}$; rather than uniformly over $Y$. 
In the same way, we can also generalise the given `global' versions of (A3) and (A5) to local versions. This will, for instance, be required in the proof of \Cref{Prop_infinite_dim_Gauss}.
However, we imagine that in most practical cases the global versions of (A3) and (A5) are not too restrictive. Hence, we for the sake of simplicity, we prefer those.
\end{remark}

\section{The additive Gaussian noise case} \label{Subsec_Gaussiancase}

In practice, the data space is typically finite dimensional and a popular modelling assumption for measurement error is additive non-degenerate Gaussian noise.
In this case, one can verify  Assumptions (A1)-(A4) --  independently of prior $\muprio$ and forward response operator $\mathcal{G} \in \mathbf{M} = \{f : X \rightarrow Y \text{ measurable}\}$. 
Hence, this very popular setting leads to a  well-posed Bayesian inverse problem in the weak topology, the Hellinger distance, and the total variation distance. 
If the prior has a finite $p$-th moment, we additionally obtain Wasserstein$(p)$ well-posedness.

\begin{corollary} \label{Cor_Gaussian_case}
Let $Y := \mathbb{R}^k$ and $\Gamma \in \mathbb{R}^{k \times k}$ be symmetric positive definite. Let $\mathcal{G} \in \mathbf{M}$ be a measurable function.  
A Bayesian inverse problem with additive non-degenerate Gaussian noise $\eta \sim \mathrm{N}(0, \Gamma)$ is given by the following likelihood:
$$L(y^\dagger|\theta) = \det(2\pi\Gamma)^{-1/2} \exp\left(-\frac{1}{2}\|\Gamma^{-1/2}(\mathcal{G}(\theta) - y^\dagger)\|_Y^2\right).$$
Then, the \cref{EQ_BIP} corresponding to likelihood $L$ and 
\begin{itemize}
\item[(a)]  any prior $\muprio \in \Prob(X)$ is whtv  well-posed;
\item[(b)] any prior $\muprio \in \Prob_p(X)$ is Wasserstein$(p)$ well-posed, where $p \in [1,\infty)$ and $X$ is a normed space.
\end{itemize}
\end{corollary}

\begin{remark} \label{Rem_Always_Wellposed_Always_Illposed}
Let $X$ contain at least two elements.
The non-Bayesian inverse problem \cref{EQ_IP} corresponding to the additive Gaussian noise setting in \Cref{Cor_Gaussian_case} is ill-posed.
We have shown this in \Cref{prop_ill-posed_IP}. 
Hence, in case of Gaussian noise, the Bayesian approach using any prior measure always gives a whtv well-posed Bayesian inverse problem, in contrast to the always ill-posed \cref{EQ_IP}.
\end{remark}

The fact that we can show well-posedness under \emph{any} forward response operator and \emph{any} prior measure in $\Prob(X)$ or $\Prob_p(X)$  has relatively strong implications for practical problems.
We now comment on the deterministic discretisation of posterior measures, hierarchical models, and Bayesian model selection.

\subsection{Deterministic discretisation}
Bayesian inverse problems can be discretised with deterministic quadrature rules; such are quasi-Monte Carlo \cite{Dick2019}, sparse grids \cite{Schwab2012}, or Gaussian quadrature. 
Those are then used to approximate the model evidence and to integrate with respect to the posterior.
Deterministic quadrature rules often behave like discrete approximations of the prior measure. 
If this discrete approximation is a probability measure supported on a finite set, we can apply \Cref{Cor_Gaussian_case} and show that the \cref{EQ_BIP} based on the discretised prior is whtv and Wasserstein$(p)$ well-posed for any $p \in [1,\infty)$.
\subsection{Hierarchical prior}
\emph{Hierachical prior measures} are used to construct more complex and flexible prior models. In Bayesian inverse problems, such are discussed in \cite{Dunlop2018,Dunlop2017,LEU19}. The basic idea is to employ a prior measure depending on a so-called hyperparameter. 
This hyperparameter has itself a prior distribution, which (typically) leads to a more complex total prior measure.
This can be continued recursively down to $K$ layers:
$$\muprio = \int_{X_K} \cdots \int_{X_1} \muprio^0( \cdot |\theta_1)\muprio^1(\mathrm{d}\theta_1|\theta_2)\ldots\muprio^K(\mathrm{d}\theta_K).$$
Here, $X_1,\ldots,X_K$ are measurable subsets of Radon spaces, $X_0 := X$, and $$\muprio^{k-1}: X_{k} \times \mathcal{B}X_{k-1} \rightarrow [0,1]$$
is a Markov kernels from $(X_k, \mathcal{B}X_{k})$ to $(X_{k-1}, \mathcal{B}X_{k-1})$, for $k \in \{1,\ldots,K\}$.
Note that hierarchical measures are in a way the probabilistic version of a deep model in machine learning --  such as a deep neural network. 
In a deep neural network, we also add layers to allow for more flexibility in function approximations.

The likelihood still depends only on $\theta$ but not on the deeper layers $\theta_1,\ldots,\theta_K$. The \cref{EQ_BIP} of determining the posterior measure $\mathbb{P}(\theta \in \cdot | y = y^\dagger)$ of the outer layer is whtv well-posed. This is a direct implication of \Cref{Cor_Gaussian_case}. 
Moreover, finding the posterior measure of all layers $\mathbb{P}((\theta,\theta_1,\ldots,\theta_K) \in \cdot | y = y^\dagger)$ is whtv well-posed, too. 
This can be seen by extending the parameter space to $X \times X_1 \times \cdots \times X_K$ 
to all layers ($\theta_k$ lives in $X_k$, $k=1,\ldots,K$) and applying \Cref{Cor_Gaussian_case} to the extended parameter space.

\subsection{Model selection}
In \emph{Bayesian model selection}, not only a model parameter shall be identified, but also the correct model in a collection of possible models. For instance, Lima et al. \cite{Lima2016} have applied Bayesian model selection to identify the correct model to represent a particular tumour. We briefly comment on a special case of Bayesian model selection.
Let $L(\cdot|\theta, \mathcal{G})$ be the likelihood in \Cref{Cor_Gaussian_case} where we now also note the dependence on the forward response operator $\mathcal{G}$.
Moreover, let $\mathbf{M}' \subseteq \mathbf{M}$ be a finite collection of forward response operators that we want to identify the \emph{correct} one from.
We now define a prior measure $\muprio'$ on $\mathbf{M}'$ which determines our a priori knowledge about the model choice.
The posterior measure of the model selection problem on $(X \times \mathbf{M}', \mathcal{B}X \otimes 2^{\mathbf{M}'})$ is given by
\begin{equation*}
\mupost^{\dagger, \rm ms} = \mathbb{P}\left((\theta, \mathcal{G}^*) \in \cdot | \mathcal{G}^*(\theta)+\eta = y^\dagger\right),
\end{equation*}
where $\mathcal{G}^*: \Omega \rightarrow \mathbf{M}'$ is the random variable representing the model; it satisfies $\mathcal{G}^* \sim \muprio'$.
The posterior can be computed using a generalisation of Bayes' Theorem
\begin{equation*}
\mupost^{\dagger, \rm ms}(A \times B) = \frac{ \sum_{\mathcal{G} \in B}\int_A L(y^\dagger|\theta, \mathcal{G}) \muprio'(\{\mathcal{G}\}) \mathrm{d}\muprio(\theta)}{\sum_{\mathcal{G}' \in \mathbf{M}'}\int_X L(y^\dagger|\theta', \mathcal{G}') \muprio'(\{\mathcal{G}'\}) \mathrm{d}\muprio(\theta')} \quad (A \in \mathcal{B}X, B \in 2^{\mathbf{M}'}).
\end{equation*}

This identity is indeed correct: We just apply \Cref{thm:bayes} on the parameter space $X \times \mathbf{M}'$ with prior measure $\muprio \otimes \muprio'$ and likelihood $L(y^\dagger|\cdot, \cdot): X \times \mathbf{M}' \rightarrow [0,\infty)$.
In the setting of \Cref{Cor_Gaussian_case}, the \cref{EQ_BIP} of identifying model and parameter is whtv well-posed.

\subsection{Generalisations}
We have discussed finite-dimen\-sional data and additive non-de\-ge\-nerate Gaussian noise. These results cannot trivially be expanded to the degenerate Gaussian noise case: Degenerate Gaussian likelihoods do not satisfy (A1) and can lead to degenerate posterior measures; we have discussed those in \Cref{SubSec_NoBayes}.

The infinite-dimensional data with additive Gaussian noise requires a likelihood definition  via the Cameron--Martin Theorem.
For a discussion of infinite-dimensional data spaces, we refer to \cite[Remark 3.8]{Stuart2010} for compact covariance operators and  \cite[\S 2.1]{Kahle2019} specifically for Gaussian white noise generalised random fields. 
Generalising the result from \cite{Kahle2019}, we can say the following:
\begin{corollary} \label{Prop_infinite_dim_Gauss}
Let $(Y', \langle \cdot, \cdot \rangle_{Y'})$ be a separable Hilbert space and $\Gamma: Y' \rightarrow Y'$ be a covariance operator; i.e. it is self-adjoint, positive definite, and trace-class.
We assume that $Y$ is the Cameron--Martin space of $\mathrm{N}(0, \Gamma) \in \Prob(Y')$, i.e. $$(Y, \langle \cdot, \cdot \rangle_Y) = (\mathrm{img}(\Gamma^{1/2},Y'), \langle \Gamma^{-1/2} \cdot,\Gamma^{-1/2} \cdot \rangle_{Y'}),$$
where the inverse square-root $\Gamma^{-1/2}$ is well-defined.
Moreover, let $\mathcal{G}: X \rightarrow Y$ be a measurable function. 
Then, the inverse problem
$$
\mathcal{G}(\theta^\dagger) + \eta = y^\dagger \qquad (\eta \sim \mathrm{N}(0, \Gamma))
$$
can be represented by the likelihood
$$
L(y^\dagger|\theta) = \exp\left(\langle \mathcal{G}(\theta), y^\dagger\rangle_Y - \frac{1}{2}\|\mathcal{G}(\theta)\|_Y^2 \right).
$$
If in addition, $\mathcal{G}:X \rightarrow Y$ is bounded, the \cref{EQ_BIP} corresponding to likelihood $L$ and 
\begin{itemize}
\item[(a)]  any prior $\muprio \in \Prob(X)$ is whtv  well-posed;
\item[(b)] any prior $\muprio \in \Prob_p(X)$ is Wasserstein$(p)$ well-posed, where $p \in [1,\infty)$ and $X$ is a normed space.
\end{itemize}
\end{corollary}

Note that we here require $\mathcal{G}$ to be bounded in $X$. Hence, while allowing for infinite-dimensional data spaces, we now have conditions on the forward response operator and also on the covariance operator. Thus, this result is not as generally applicable as \Cref{Cor_Gaussian_case}.

Generalisations to non-Gaussian noise models are not as difficult. In the proof of \Cref{Cor_Gaussian_case}, we have only used that the probability density function of the noise is strictly positive, continuous in its argument, and bounded by a constant.
This however is also satisfied, when the noise is additive, non-degenerate  and follows, e.g., the Cauchy distribution, the $t$-distribution, or the Laplace distribution.

\section{Numerical illustrations} \label{Sec_Num_illust}

We illustrate some of the results shown in the sections before with numerical examples.
First, we consider some simple one-dimensional examples complementing the examples we have considered throughout the article.
Those include Bayesian inverse problems with likelihoods that are discontinuous in parameter or data.
Second, we consider an inverse problem that is high-dimensional in terms of data and parameters.
The high-dimensional inverse problem is concerned with the reconstruction of an image by Gaussian process regression.
\subsection{Discontinuities in the likelihood}
\begin{figure}
\centering
\includegraphics[scale=0.8]{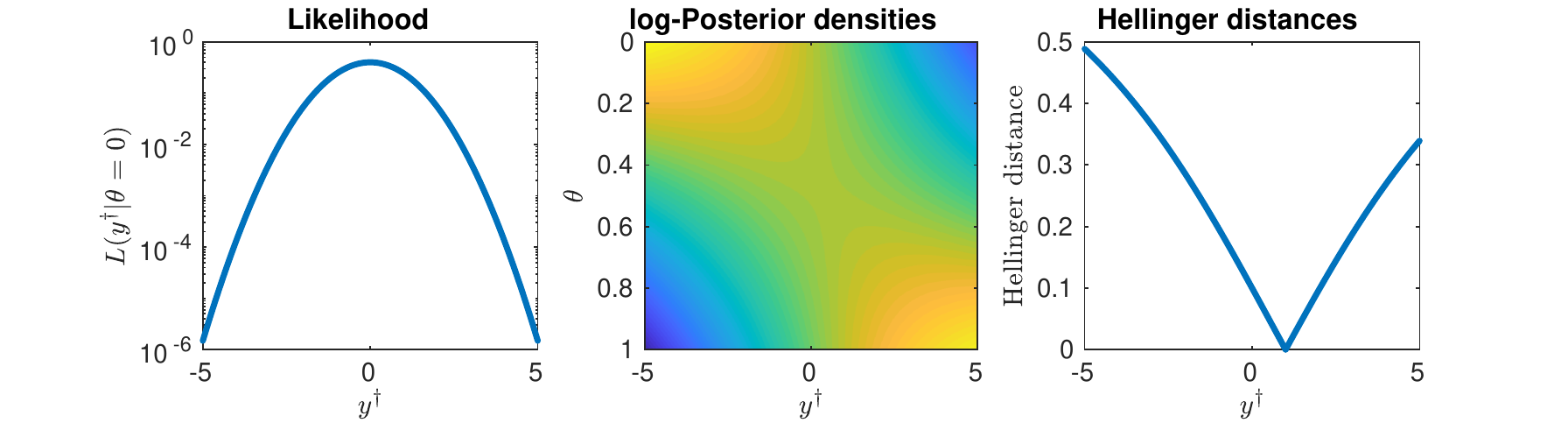}
\includegraphics[scale=0.8]{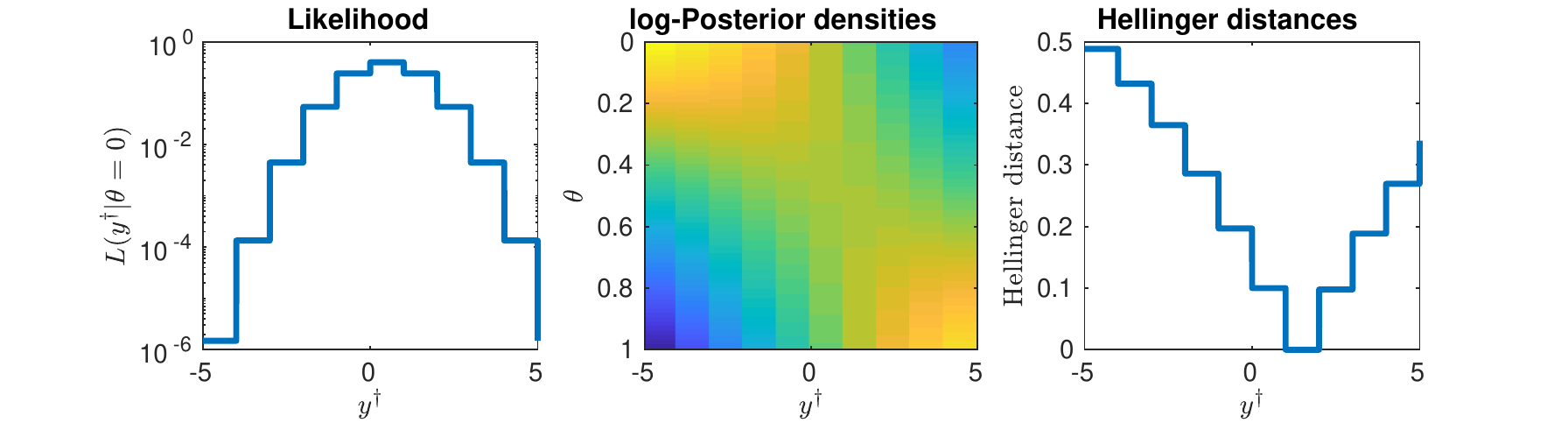}
\caption{Top row:  \Cref{Ex_nonctscts_Like}(a), bottom row:  \Cref{Ex_nonctscts_Like}(b). Left: Likelihood at $\theta = 0$. Centre: Log-posterior densities corresponding to the Bayesian inference problems. The colormaps show a descent in posterior density, when going from yellow (high) to dark blue (low). Right: Hellinger distance between the posterior $\mupost^\ddagger$ with $y^\ddagger = 1$ and posterior $\mupost^\dagger$ with $y^\dagger$ varying between $-5$ and $5$.}
\label{Fig_likeli_non_cts}
\end{figure}

In previous works, Lipschitz continuity of the log-likelihood in the data and (at least) continuity in the parameter has been assumed, see \cite{Stuart2010}. 
In this article, we prove results that do not require continuity in the parameter, however, we still require continuity in the data.
We now illustrate these results with simple numerical experiments.
Indeed, we show that Assumption (A4) is crucial, by comparing \cref{EQ_BIP} posteriors with likelihoods that are continuous and  discontinuous in the data.
\begin{example}[Continuity of $y \mapsto L(y| \cdot)$] \label{Ex_nonctscts_Like}
We define data and parameter space by $Y := \mathbb{R}$ and $X:=[0,1]$. We consider the \cref{EQ_BIP}s with prior measure $\muprio := \mathrm{Unif}(0,1)$ on $X$ and one of the following likelihoods
\begin{itemize}
\item[(a)] $L(y^\dagger|\theta) = (2 \pi)^{-1/2} \exp(-\frac12 \|y^\dagger - \theta\|_Y^2)$, 
\item[(b)] $L(y^\dagger|\theta) = (2 \pi)^{-1/2} \exp(-\frac12 \|\lfloor y^\dagger \rfloor - \theta\|_Y^2)$.
\end{itemize}
\end{example}
We solve the inverse problems in \Cref{Ex_nonctscts_Like} with numerical quadrature. 
In particular, we compute the model evidences for a $y^\dagger \in \{-5,-4.999,-4.998,\ldots5\}$ and the Hellinger distances between $\mupost^\dagger$ and $\mupost^\ddagger$, where $y^\ddagger = 1$.
In \Cref{Fig_likeli_non_cts}, we plot the likelihood functions at $\theta = 0$, the logarithms of the posterior densities, and the Hellinger distances.
The top row in the figure refers to (a), the bottom row refers to (b).
In the continuous setting (a), we see continuity with respect to $y^\dagger$ in all images. 
Indeed the \cref{EQ_BIP} in (a) fulfills (A1)-(A4). 
The inverse problem in (b) satisfies (A1)-(A3), but not (A4).
Also, we see discontinuities with respect to the data in all figures referring to (b).
Especially, the figure of the Hellinger distances is discontinuous which leads to the conclusion that this inverse problem is not  well-posed.
Hence, (A4) is indeed crucial to obtain  well-posedness of a Bayesian inverse problem.

\begin{remark}
A likelihood as in \Cref{Ex_nonctscts_Like}(b) can arise, when considering cumulative or categorial data, rather than real-valued continuous data as in (a).
Categorial data arises in classification problems.
\end{remark}

While continuity in the data is important, we now illustrate that continuity in the forward response operator is not necessary to obtain continuity in the data to posterior map.
We give an example that can be understood as learning the bias in a single layer neural network.
\begin{example}[Continuity in $\theta \mapsto L(\cdot|\theta)$] \label{Ex_nonctscts_forward }
We define data and parameter space by $Y := \mathbb{R}$ and $X:=[0,1]$.
Let $w \in [1,\infty]$ be a known weight parameter.
We define the forward response operator with weight $w$ by
$$\mathcal{G}_w: X \rightarrow Y, \qquad \theta \mapsto \frac{1}{1+\exp(-w(0.5-\theta))}.$$
If $w < \infty$, the forward response operator resembles a single layer neural network with sigmoid activation function evaluated at $0.5$. This neural network has known weight $w$ and uncertain bias $\theta$.
Moreover, note that in the limiting setting $w= \infty$, the sigmoid function is there replaced by the heaviside function with step at $\theta$, evaluated also at $x = 0.5.$
\begin{equation}
\mathcal{G}_\infty: X \rightarrow Y, \qquad \theta \mapsto \begin{cases} 1, &\text{ if } 0.5 \geq \theta, \\ 0, &\text{ otherwise}. \end{cases}
\end{equation}
We consider the \cref{EQ_BIP} of estimating the true bias $\theta^\dagger$, given an observation $y^\dagger_w := \mathcal{G}_w(\theta^\dagger) + \eta^\dagger$.
Here, we consider the noise $\eta^\dagger$ to be a realisation of $\eta \sim \mathrm{N}(0,1^2)$.
Moreover, we assume that the parameter $\theta \sim \muprio =  \mathrm{Unif}(0,1)$ follows a uniform prior distribution.
\end{example}
\begin{figure}
\includegraphics[scale=0.7]{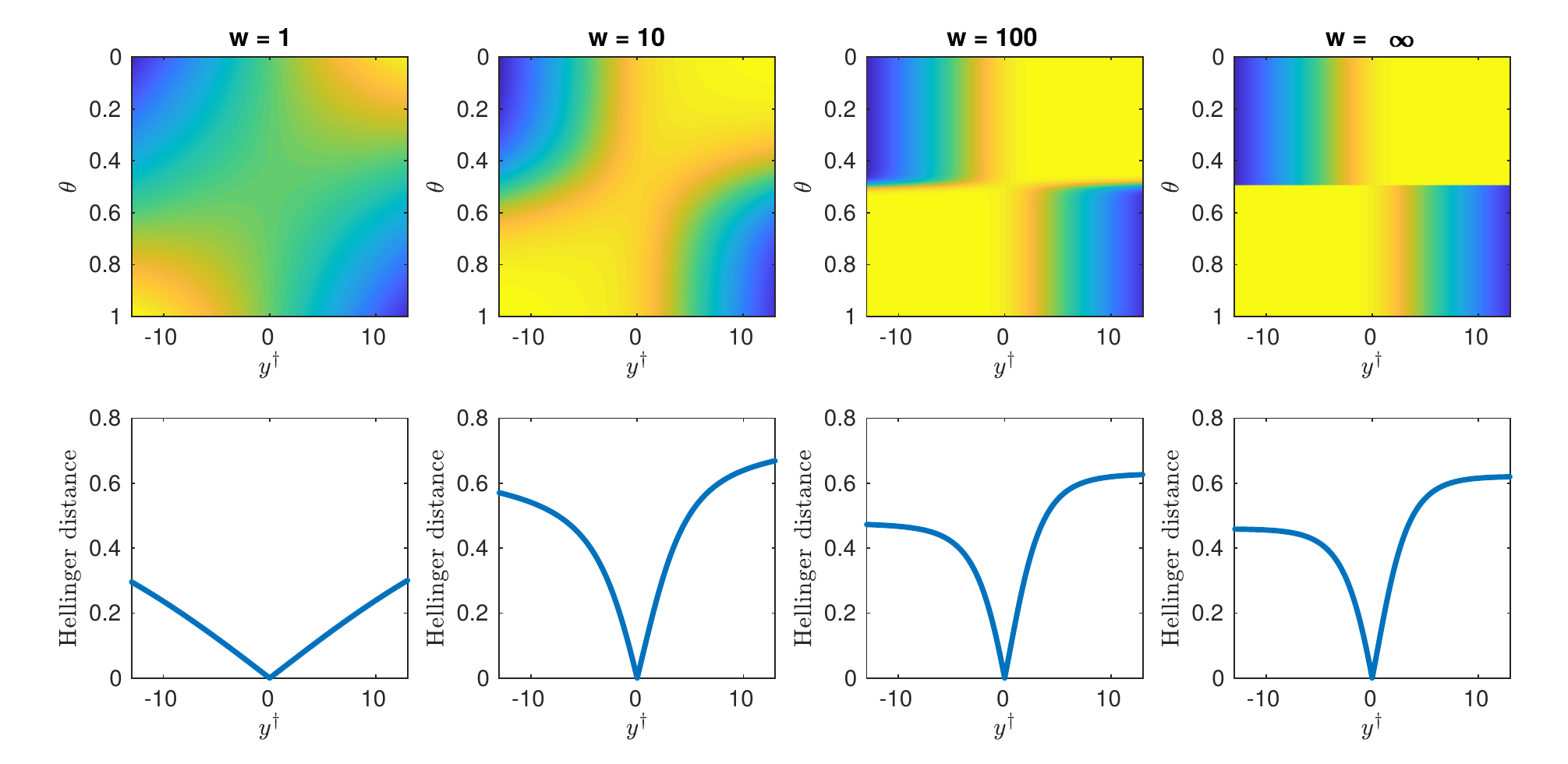}
\caption{From left to right: \Cref{Ex_nonctscts_forward } given $w=1,10,100,\infty$. Top row: Log-posterior densities corresponding to the Bayesian inference problems. The colormaps show a descent in posterior density, when going from yellow (high) to dark blue (low). Bottom row: Hellinger distance between the posterior $\mupost^\ddagger$ with $y^\ddagger = 1$ and posterior $\mupost^\dagger$ with $y^\dagger$ varying between $-13$ and $13$.}
\label{Fig_forward_non_cts}
\end{figure}
We solve the \cref{EQ_BIP}s in \Cref{Ex_nonctscts_forward } with weights $w = 1,10,100,\infty$ again with numerical quadrature for $y^\dagger \in \{-13,-12.99,-12.98,\ldots,13\}$.
We compute the Hellinger distance between $\mupost^\dagger$ and $\mupost^\ddagger$, where $y^\ddagger=0$.
We plot the logarithms of the posterior densities of obtained in \Cref{Ex_nonctscts_forward } in \Cref{Fig_forward_non_cts}, along with the Hellinger distances.
We observe that all of the posteriors are continuous with respect to the data. 
This includes the posterior that is based on the discontinuous forward response operator $\mathcal{G}_\infty$.
It is discontinuous in the parameter, but continuous in the data. 
The \cref{EQ_BIP} considered here satisfy again (A1)-(A4).
Hence, also these numerical experiments verify the statement of \Cref{Lemma_main}.

\begin{remark}
In \emph{deep learning}, sigmoid functions $\mathcal{G}_w$ $(w < \infty)$ are considered as smooth approximations to the heaviside function $\mathcal{G}_\infty$, which shall be used as an activation function. 
The smooth sigmoid functions allow to train the deep neural network with a gradient based  optimisation algorithm.
When training the neural network with a Bayesian approach, rather than an optimisation approach, we see that we can use heaviside functions in place of smooth approximations and obtain a well-posed Bayesian inverse problem.
\end{remark}

\subsection{A high-dimensional inverse problem} \label{Subsec_Cat}
We now consider an inverse problem that is high-dimensional in parameter and data space. 
In particular, we observe single, noisy pixels of a  grayscale photograph. 
The inverse problem consists in the reconstruction of the image, for which we use Gaussian process regression. 
We then perturb the data by adding white noise to the image and investigate changes in the posterior, as we rescale the noise.

\begin{figure}
\centering
\includegraphics[scale=0.65]{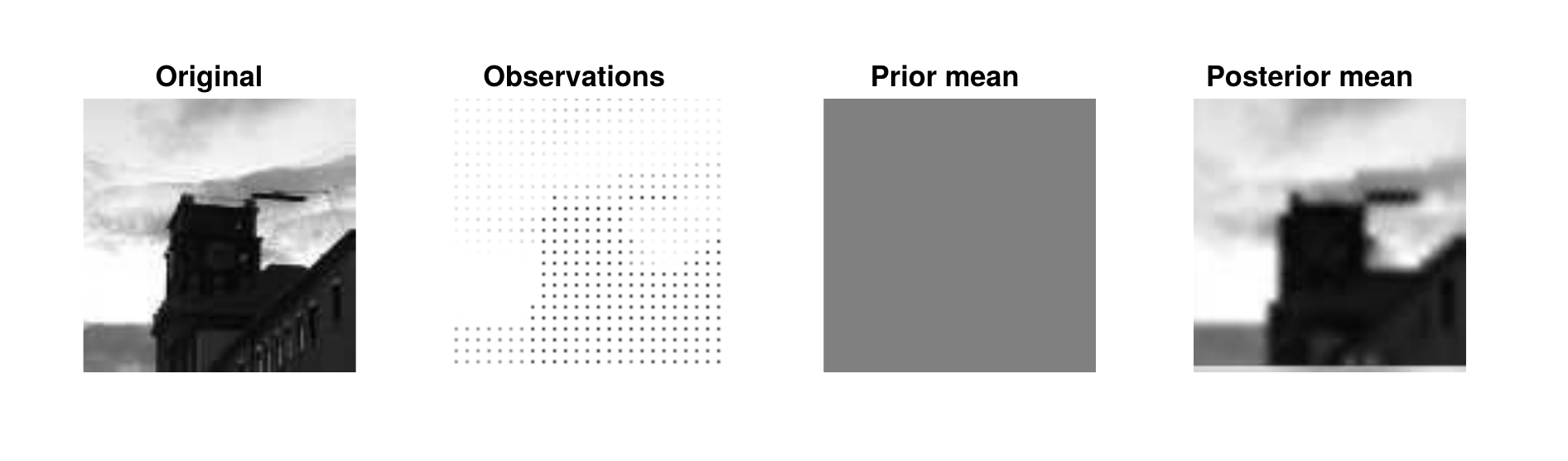}
\caption{Reconstruction of an image with Gaussian process regression. From left to right: original image, observational data (white parts are unobserved), prior mean, and posterior mean}
\label{Fig_Klaus_reconstruction}
\end{figure}

\begin{figure}
\centering
\includegraphics[scale=0.65]{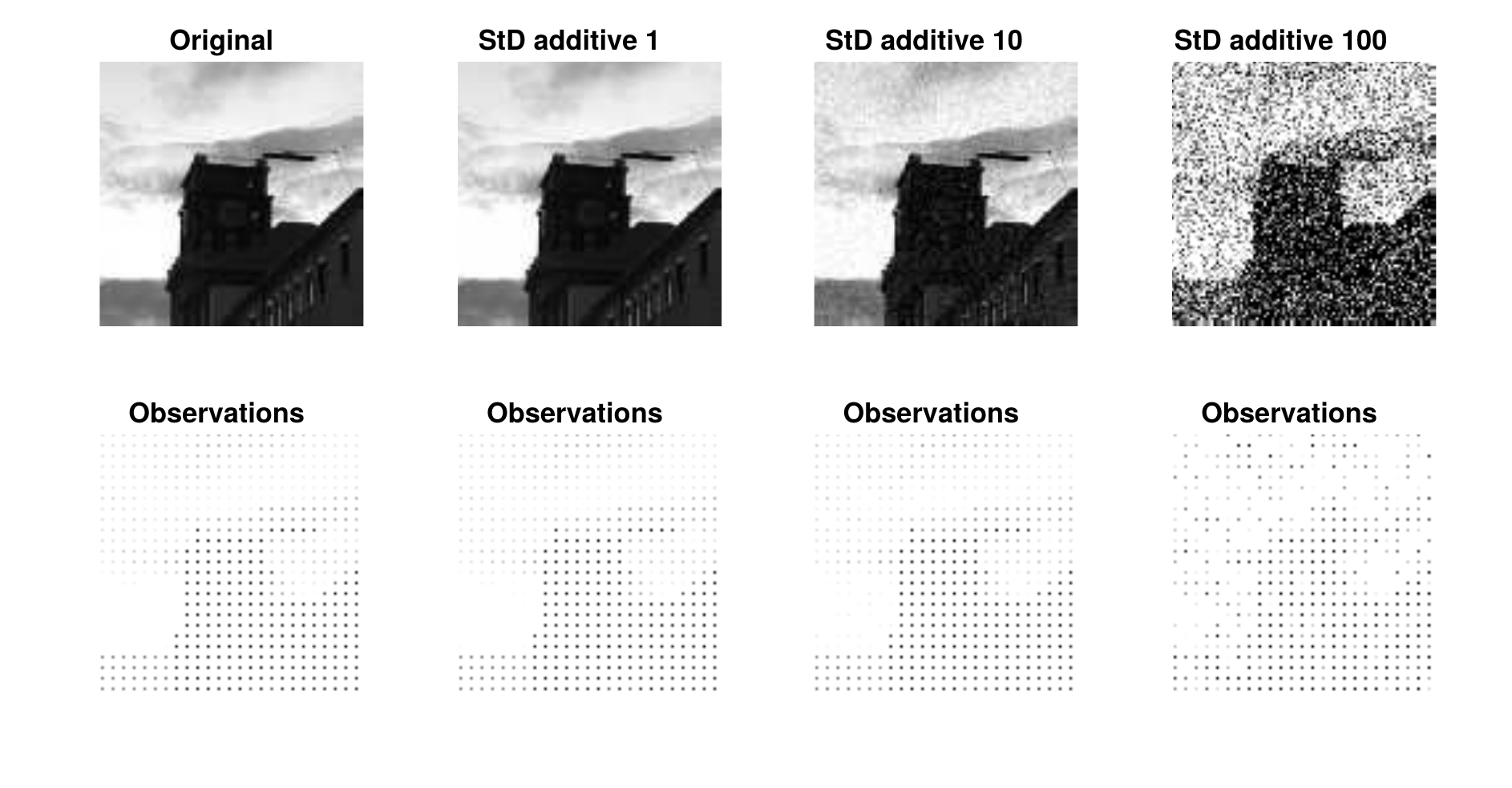}
\caption{Top row: Original image and images perturbed with scaled white noise, given $\sigma \in \{1,10,100\}$. Bottom row: Observations obtained from the perturbed image.}
\label{Fig_Klaus_rotation}
\end{figure}

\begin{example} \label{Ex_cat_image}
Let the parameter space $X := \mathbb{R}^{100 \times 100}$ contain grayscale images made up of $100 \times 100$ pixels. The data space $Y := \mathbb{R}^{25 \times 25}$ consists of $25 \times 25$ pixels that are observed in a single picture. 
Returning those  $25 \times 25$ pixels from a $100 \times 100$ pixels image is modelled by the function $\mathcal{G}:X \rightarrow Y$.
Let $\theta^\dagger \in X$ be a full image. 
Given $$y^\dagger = \mathcal{G}(\theta^\dagger)+\eta,$$
we shall recover the full image $\theta^\dagger$. 
Here, $\eta \sim \mathrm{N}(0, 5^2 I)$ is normally distributed noise, with a noise level of about $5/\max(y) = 2\%$.
We assume a Gaussian prior on $X$:
$$\muprio = \mathrm{N}\left(\begin{pmatrix}
128 &\cdots&128 \\
\vdots & \ddots & \vdots \\
128 &\cdots &128
\end{pmatrix}, C_0 \right),$$
where $C_0 \in \mathbb{R}^{100^{\times 4}}$ is a covariance tensor assigning the following covariances:
$$\mathrm{Cov}(\theta_{i,j},\theta_{\ell,k}) = 10000\cdot\exp\left(-\frac{\sqrt{(i-\ell)^2 +(j-k)^2}}{15}\right).$$
Note that this is essentially an adaptation of an exponential covariance kernel for a  Gaussian process in 2D space.
\end{example}
The Bayesian inverse problem in \Cref{Ex_cat_image} can be solved analytically, since $\mathcal{G}$ is linear, and prior and noise are Gaussian.
We obtain the posterior measure by Gaussian process regression.
In \Cref{Fig_Klaus_reconstruction}, we present the original image, observations, prior mean image and posterior mean image.
The reconstruction is rather coarse, which is not surprising given that we observe only $6.25 \cdot 10^2$ of $10^4$ pixels of the image.  
We now investigate how the posterior measure changes under marginal changes in the data.
To do so, we perturb the image additively with  scaled white noise. In particular, we add $\mathrm{N}(0, \sigma^2)$-distributed, independent random variables to each pixel.
In \Cref{Fig_Klaus_rotation}, we show images and associated observations, where the standard deviations (StD) of the noise is $\sigma \in \{1, 10, 100\}$.

\begin{figure}
\centering
\includegraphics[scale=0.55]{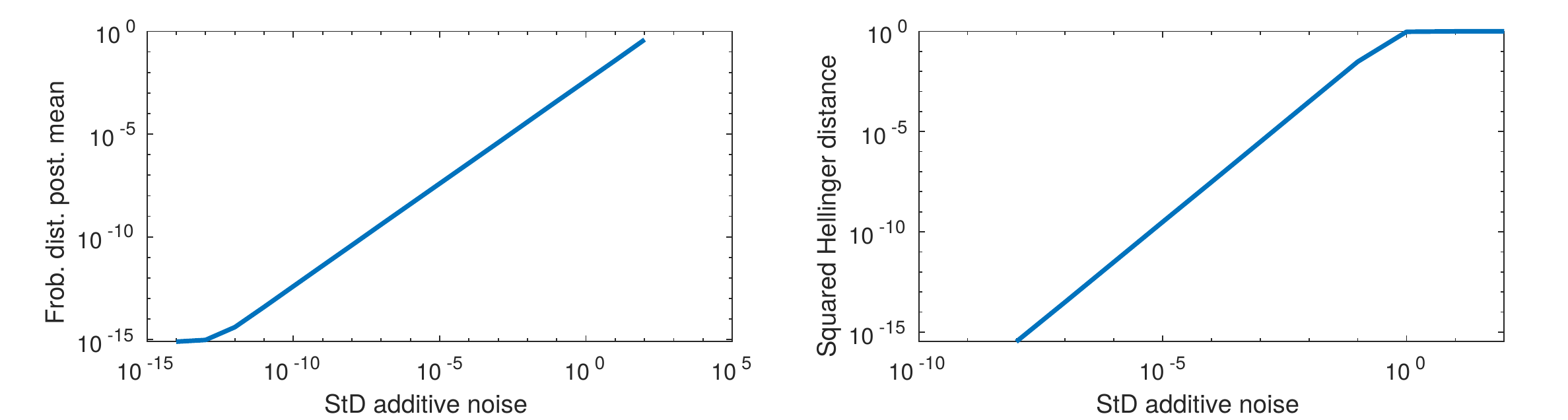}
\caption{Mean relative Frobenius distances and mean squared Hellinger distances computed between the posterior $\mupost^\dagger$ and posteriors $\mupost^\ddagger$, in which the underlying image was perturbed with white noise that has been scaled by StD $\sigma = 0, 10^{-17}, 10^{-19},\ldots,10^{2}$. `Mean' refers to the fact that the perturbations are random and the distances have been computed for 20 random perturbations and then averaged.
When approaching $|y^\ddagger - y^\dagger| \rightarrow 0$, the distances go to $0$. The left-out $x$-values have distance zero numerically. }
\label{Fig_Klaus_continuity}
\end{figure}

Using Gaussian process regression, we compute the posteriors after perturbing the images with scaled white noise given $\sigma \in \{10^{-17}, 10^{-16},\ldots,10^{2}\}$.
Between the original posterior with no perturbation in the data and all others, we compute the Hellinger distance and the relative Frobenius distance between the (matrix-valued) posterior means
$$\text{Relative Frobenius distance } = \frac{\Big\|\int \theta \mathrm{d}\mupost^\ddagger(\theta) - \int \theta \mathrm{d}\mupost^\dagger(\theta)\Big\|_{\rm F}}{\Big\|\int \theta \mathrm{d}\mupost^\ddagger(\theta)\Big\|_{\rm F}},
$$
where $\mupost^\dagger$ (resp. $\mupost^\ddagger$) is the posterior referring to the perturbed data  $y^\dagger$ (resp. non-perturbed data $y^\ddagger$).
Since the perturbation is random, we perform this process 20 times and compute the mean over these distances.
The standard deviation in these metrics is negligibly small.
We plot the results in \Cref{Fig_Klaus_continuity}, where we see indeed continuity reducing the error standard deviation in the data.
In light of \Cref{Lemma_main} and \Cref{Cor_Gaussian_case}, this is what we expect: First, note that the Bayesian inverse problem falls in the category \emph{additive finite-dimensional Gaussian noise} and is therefore  well-posed.
Hence, also in this high-dimensional setting, we are able to verify our analytical results concerning well-posedness.

\section{Conclusions}\label{Sec_Conclusions}
\begin{figure}[hptb]
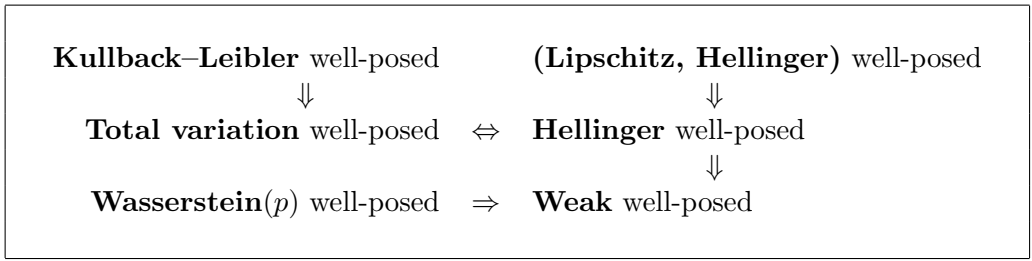

\begin{center}
\begin{tabular}{|rrcll|}\hline
\quad & &\qquad & &\quad \\
\quad & \textbf{Kullback--Leibler} well-posed & & \textbf{(Lipschitz, Hellinger)} well-posed &\quad \\
\quad &   $\Downarrow$ \qquad \qquad \qquad  && \qquad  \qquad  \qquad  $\Downarrow$  &\quad \\
\quad & \textbf{Total variation} well-posed    &    
$\Leftrightarrow$  &                     
   \textbf{Hellinger} well-posed \qquad   &\quad  \\  
\quad & & & \qquad  \qquad  \qquad  $\Downarrow$               &\quad         \\
\quad &\textbf{Wasserstein}$(p)$ well-posed & $\Rightarrow$ &  \textbf{Weak} well-posed   
&\quad\\

\quad & &\qquad & &\quad \\
\hline
\end{tabular}
\end{center}
\caption{Relations between concepts of  well-posedness. Here, $A \Rightarrow B$ means that an \cref{EQ_BIP} being $A$-well-posed implies that it is also $B$-well-posed.}
\label{Figure_Implicationscheme}
\end{figure}

In this work, we introduce and advocate a new concept of well-posedness of Bayesian inverse problems.
We weaken the stability condition by considering continuity instead of Lipschitz continuity of the data-to-posterior map. 
On the other hand, we make the stability condition somewhat stronger by allowing to adapt the metric on the space of probability measures to the particular situation.
From this discussion arise various notions of well-posedness, which and whose relations we summarise in \Cref{Figure_Implicationscheme}.

Importantly, we show that, given our concept, a huge class of practically relevant Bayesian inverse problems is well-posed or can easily be shown to be well-posed.
Hence, we give a general justification for the Bayesian approach to inverse problems for a huge number of practical situations: here, the Bayesian inverse problem will have a unique solution and this unique solution will be robust with respect to marginal changes in the data. Such inverse problems appear e.g. in engineering, machine learning, and imaging.

Concerning directions for future research, we denote the following. So far, we have mostly neglected the degenerate Bayesian inverse problems which we discussed in \Cref{SubSec_NoBayes}. Such problems appear in Bayesian probabilistic numerics, and other settings where noise-free data is considered.
This may also include the Bayesian formulation of machine learning problems with discrete loss models, like \emph{0-1-loss}, or Bayesian formulations of classification problems, see \cite{Bertozzi2018}.

\appendix

\section{Conditional probability} \label{Subsec_Appedix_CondProb} In this appendix, we briefly summarise some results concerning conditional probabilities. 
Let $X, Y$ be given as in \Cref{SubSec_IPs}. Moreover, let $\Omega := X \times Y$ and $\theta : \Omega \rightarrow X, y: \Omega \rightarrow Y$ be random variables.
\begin{theorem} \label{Thm_RegCondProb}
A Markov kernel $M: Y \times  \mathcal{B}X \rightarrow [0,1]$ exists, such that
$$
\mathbb{P}(\{\theta \in A\} \cap \{y \in C\}) = \int_C M(y^\dagger, A) \mathbb{P}(y \in \mathrm{d}y^\dagger) \qquad (A \in \mathcal{B}X, C \in \mathcal{B}Y)
$$
Moreover, $M$ is $\mathbb{P}(y \in \cdot)$-a.s. unique.
\end{theorem} 

Let $y^\dagger \in Y$.
The probability measure $M(y^\dagger, \cdot)$ in \Cref{Thm_RegCondProb} is the \emph{(regular) conditional probability distribution} of $\theta$ given that $y = y^\dagger$. We denote it by $\mathbb{P}(\theta \in \cdot | y = y^\dagger)$. Note that the conditional probability is only unique for a.e. $y^\dagger \in Y$.
This definition, as well as \Cref{Thm_RegCondProb}, is non-constructive. 
However, if we can represent the joint distribution $\mathbb{P}((\theta, y) \in \cdot)$ by a probability density function, we can compute the density of the conditional probability distribution.
First, consider the following lemma concerning joint and marginal probability density functions.
\begin{lemma} \label{Lemma_marginalpdf}
Let $\nu_X$ and $\nu_Y$ be $\sigma$-finite measures on $(X, \mathcal{B}X)$ and $(Y, \mathcal{B}Y)$ and $$\mathbb{P}((\theta, y) \in \cdot) \ll \nu_X \otimes \nu_Y, \text{with } f := \frac{\mathrm{d}\mathbb{P}((\theta, y) \in \cdot)}{\mathrm{d}\nu_X \otimes \nu_Y} \qquad (\nu_X \times \nu_Y\text{-a.e.}).$$ 
Then, 
$\mathbb{P}(\theta \in \cdot) \ll \nu_X,$ with ${\mathrm{d}\mathbb{P}(\theta \in \cdot)}/{\mathrm{d}\nu_X} = \int_X f(\cdot, y^\dagger) \nu_Y(\mathrm{d}y^\dagger), \nu_X\text{-a.e.}$,  and  $\mathbb{P}(y \in \cdot) \ll \nu_Y$ with ${\mathrm{d}\mathbb{P}(y \in \cdot)}/{\mathrm{d}\nu_Y} = \int_Y f(\theta^\dagger, \cdot ) \nu_X(\mathrm{d}\theta^\dagger), \nu_Y\text{-a.e.}$
\end{lemma}

Next, we move on to the construction of the conditional density.

\begin{lemma}\label{Lemma_condpdf}
Let $\nu_X, \nu_Y$, and $f$ be given as in \Cref{Lemma_marginalpdf}. Then, for $\theta^\dagger \in X$ $(\nu_X\text{-a.e.})$ and $y^\dagger \in Y$ $(\nu_Y\text{-a.e.})$, we have
$$\frac{\mathrm{d}\mathbb{P}(\theta \in \cdot | y = y^\dagger)}{\mathrm{d}\nu_X}(\theta^\dagger) = \begin{cases} \frac{f(\theta^\dagger, y^\dagger)}{g(y^\dagger)}, &\text{if } g(y^\dagger) > 0, \\ 0, &\text{otherwise},\end{cases}$$
where $g(y^\dagger) := \int_X f(\theta^\ddagger, y^\dagger) \nu_X(\mathrm{d}\theta^\ddagger)$ is the $\nu_Y$-probability density function of $\mathbb{P}(y \in \cdot)$.
\end{lemma}

This result is fundamental to prove Bayes' Theorem; see \Cref{thm:bayes}.

\section{Proofs} \label{Subsec_Appedix_Proofs} In this appendix, we present rigorous proofs of all the theorems, propositions, lemmata, and corollaries stated in this article.
\begin{proof}[Proof of \Cref{prop_ill-posed_IP}]
Note that the support of $\munoise$ is all of $Y$. 
Hence, the noise $\eta^\dagger$ can be any value in $Y$ and we need to solve the equation
\begin{equation}\label{Eq_Proof_Illposed}
y^\dagger = \mathcal{G}(\theta^\dagger) + \eta^\dagger
\end{equation}
with respect to both, $\theta^\dagger \in X$ and $\eta^\dagger \in Y$.
Let $\theta' \in X$.  Set $\eta' := y^\dagger- \mathcal{G}(\theta')$. Then, $(\theta', \eta')$ solves Equation \cref{Eq_Proof_Illposed} and thus the inverse problem \cref{EQ_IP}. 
Hence, each element in $X$ implies a solution.
Since $X$ contains at least two elements, the solution is not unique and, thus, \cref{EQ_IP} is ill-posed.
\end{proof}

\begin{proof}[Proof of \Cref{thm:bayes}] The following statements hold $\mathbb{P}(y \in \cdot)$-a.s. for $y^\dagger \in Y$.

We first show that $Z(y^\dagger) > 0$.
Since we assume that $L(y^\dagger|\cdot)$ is $\muprio$-a.s. strictly positive, we can write:
\begin{align}\label{EQ_Bayes_proof_model_evidence}
Z(y^\dagger) = \int L(y^\dagger|\theta) \mathrm{d}\muprio(\theta) 
= \int_{\{L(y^\dagger| \cdot ) > 0\}} L(y^\dagger|\theta)\mathrm{d}\muprio(\theta).
\end{align}
Now let $n \in \mathbb{N}$. As the integrand in \eqref{EQ_Bayes_proof_model_evidence} is positive, Chebyshev's inequality, \cite[Theorem 2.4.9]{Ash2000}, implies that
\begin{equation} \label{EQ_Markovs_ineq}
n \cdot \int_{\{L(y^\dagger| \cdot ) > 0\}} L(y^\dagger|\theta)\mathrm{d}\muprio(\theta) \geq  \muprio(L(y^\dagger| \cdot ) > n^{-1}).
\end{equation}
We aim to show that the probability on the right-hand side of this equation converges to $1$ as $n \rightarrow \infty$. 
Knowing this, we can conclude that  the right-hand side is strictly positive for all $n \geq N$, for some $N \in \mathbb{N}$.

Note that measures are continuous with respect to increasing sequences of sets. 
We define the set
$$B_n := \{L(y^\dagger| \cdot ) > n^{-1}\}$$
and observe that 
 $(B_n)_{n=1}^\infty$ is indeed an increasing sequence.
Moreover, note that 
$$B_\infty = \bigcup_{m=1}^\infty B_m = \{L(y^\dagger| \cdot ) > 0\},$$
and that $\muprio(B_\infty) = 1$.
Hence, we have $$\lim_{n \rightarrow \infty} \muprio(L(y^\dagger| \cdot ) > n^{-1}) = \muprio(L(y^\dagger| \cdot ) > 0) = 1.$$
As mentioned earlier, we now deduce that for some $\varepsilon \in (0, 1)$, there is an index $N \in \mathbb{N}$ such that
$$
|\muprio(L(y^\dagger| \cdot ) > n^{-1}) - 1| \leq \varepsilon < 1 \qquad (n \geq N)
$$
and thus $\muprio(L(y^\dagger| \cdot ) > n^{-1}) > 0$, for $n \geq N$. 
Plugged into Equation \eqref{EQ_Markovs_ineq}, this gives us $Z(y^\dagger) > 0.$ We have also $Z(y^\dagger) < \infty$, since $L(y^\dagger|\cdot) \in \mathbf{L}^1(X,\muprio)$, Thus, the posterior density  \eqref{EQ_Bayes} is well-defined.
We now apply Bayes' Theorem in the formulation of  \cite[Theorem 3.4]{Dashti2017} and obtain
$$\frac{\mathrm{d}\mupost^\dagger}{\mathrm{d}\muprio}(\theta') = \frac{L(y^\dagger|\theta')}{Z(y^\dagger)}, \qquad (\theta' \in X, \muprio\text{-a.s.}).$$
This implies 
$$\pipost^\dagger(\theta') = \frac{\mathrm{d}\mupost^\dagger}{\mathrm{d}\nu_X}(\theta')  = \frac{\mathrm{d}\mupost^\dagger}{\mathrm{d}\muprio}(\theta')  \frac{\mathrm{d}\muprio}{\mathrm{d}\nu_X}(\theta')  =  \frac{L(y^\dagger|\theta')\piprio(\theta')}{Z(y^\dagger)}, \quad (\theta' \in X, \nu_X\text{-a.s.}),$$
by application of standard results concerning Radon-Nikodym derivatives. This concludes the proof.
\end{proof}

\begin{proof}[Proof of \Cref{prop_homeomo}]
We test $\mupost^\dagger$ in \Cref{Thm_RegCondProb}. Let $\theta \sim \muprio$ and $y \sim \mu_L(\cdot|\theta)$. Then, $\mathbb{P}(y = \mathcal{G}(\theta)) = 1.$ Therefore, for $A \in \mathcal{B}X, C \in \mathcal{B}Y$, we have
\begin{align*}
\mathbb{P}(\{\theta \in A\} \cap \{y \in C\}) &=  \mathbb{P}(\{y \in \mathcal{G}(A)\} \cap \{y \in C\}) \\
&= \int_C \mathbf{1}_{\mathcal{G}(A)}(y^\dagger) \mathbb{P}(y \in \mathrm{d}y^\dagger) \\
&= \int_C \delta\left(A - \mathcal{G}^{-1}(y^\dagger)\right) \mathbb{P}(y \in \mathrm{d}y^\dagger).
\end{align*}
Note that $\mathcal{G}(A) \in \mathcal{B}Y$, since $\mathcal{G}^{-1}$ is continuous.
Hence, according to \Cref{Thm_RegCondProb}, we have $\mathbb{P}(\theta \in \cdot| \mathcal{G}(\theta) = y^\dagger) = \delta\left(\cdot - \mathcal{G}^{-1}(y^\dagger)\right)$, for   $\mathbb{P}(y \in \cdot)$-a.e. $y^\dagger \in Y$. Moreover, we have $\mathbb{P}(y \in \cdot) = \muprio(\mathcal{G} \in \cdot)$.
\end{proof}

\begin{proof}[Proof of \Cref{Thm_main}]
Hellinger well-posedness follows from \Cref{Lemma_main}. There, we show existence and uniqueness on $\PP := \Prob(X,\muprio)$. 
According to \Cref{thm:bayes}, we again obtain existence and uniqueness of the posterior measure also on $\PP := \Prob(X)$, as required for weak and total variation well-posedness. 
By \cite{Gibbs2002}, we have
$$ \dProk(\mu, \mu') \leq  \dtv(\mu, \mu') \leq \sqrt{2} \dHel(\mu, \mu') \qquad (\mu, \mu' \in \Prob(X,\muprio)).$$
Hence, $\dProk$ and $\dtv$ are coarser than $\dHel$. By \Cref{Propo_wellposedd1d2}, the \cref{EQ_BIP} is weakly and total variation well-posed.
\end{proof}

\begin{proof}[Proof of \Cref{Lemma_main}]
Note that existence and uniqueness of the measure $\mupost^\dagger$ are results of \Cref{thm:bayes} that holds, since (A1)-(A2) are satisfied.
We proceed as follows: we show that the likelihood is continuous as a function from $Y$ to $\mathbf{L}^1(X, \muprio)$ and that at the same time $y^\dagger \mapsto Z(y^\dagger)$ is continuous.
This implies that $y^\dagger \mapsto L(y^\dagger|\cdot)^{1/2} \in \mathbf{L}^2(X, \muprio)$ is continuous as well.
Then, we collect all of this information and show the continuity in the Hellinger distance, which is the desired result.

1. We now show continuity in $y^\dagger \in Y$ when integrating $L(y^\dagger|\cdot)$ with respect to $\muprio$.
This is  a standard application of \emph{Lebesgue's dominated convergence theorem} (DCT):
Let $(y_n)_{n=1}^\infty \in Y^\mathbb{N}$ be a sequence converging to $y^\dagger$, as $n\rightarrow \infty$. 
(A4) implies that $\lim_{n\rightarrow \infty}L(y_n|\cdot) = L(y^\dagger|\cdot)$ pointwise in $X$.
We obtain by the DCT
\begin{align*}
\lim_{n \rightarrow \infty} \int L(y_n|\cdot) \mathrm{d}\muprio =  \int\lim_{n \rightarrow \infty} L(y_n|\cdot) \mathrm{d}\muprio = \int L(y^\dagger|\cdot) \mathrm{d}\muprio,
\end{align*}
since the sequence $(L(y_n|\cdot))_{n=1}^\infty$ is bounded from above by $g \in \mathbf{L}^1(X, \muprio)$ and bounded from below by $0$, see (A1) and (A3).
Hence, the functions $$Y \ni y^\dagger \mapsto \int L(y^\dagger|\cdot) \mathrm{d}\muprio = Z(y^\dagger) \in \mathbb{R}, \qquad Y \ni y^\dagger \mapsto L(y^\dagger|\cdot) \in \mathbf{L}^1(X, \muprio)$$ are continuous.
Moreover, note that \Cref{thm:bayes} implies  that $Z(y^\dagger) $ is finite and strictly larger than $0$.

2. The continuity in $\mathbf{L}^1(X, \muprio)$ implies that for every $y^\dagger \in Y$, we have for $\varepsilon_1 > 0$ some $\delta_1(\varepsilon_1) > 0$, such that 
$$\|L(y^\dagger|\cdot) - L(y^\ddagger|\cdot)\|_{\mathbf{L}^1(X, \muprio)} \leq \varepsilon_{1} \quad (y^\ddagger \in Y: \|y^\dagger - y^\ddagger\|_Y \leq \delta_{1}(\varepsilon_1)).$$
Using this, we can  show that $y^\dagger \mapsto L(y^\dagger|\cdot)^{1/2}$ is continuous in $\mathbf{L}^2(X, \muprio)$.
Let $y^\dagger \in Y$ and $\varepsilon_1, \delta_1(\varepsilon_1), y^\ddagger$ be chosen as above. 
We have 
\begin{align*}
\|L(y^\dagger|\cdot)^{1/2} &- L(y^\ddagger|\cdot)^{1/2}\|_{\mathbf{L}^2(X, \muprio)}^2 \\ &= \int \big|L(y^\dagger|\cdot)^{1/2} - L(y^\ddagger|\cdot)^{1/2} \big|^2\mathrm{d}\muprio \\
&\leq \int \big|L(y^\dagger|\cdot)^{1/2} - L(y^\ddagger|\cdot)^{1/2} \big|\times\big|L(y^\dagger|\cdot)^{1/2} + L(y^\ddagger|\cdot)^{1/2} \big|\mathrm{d}\muprio \\
&= \int\big|L(y^\dagger|\cdot) - L(y^\ddagger|\cdot) \big|\mathrm{d}\muprio \leq \varepsilon_1. \\
\end{align*}
Now, we take the square-root on each side of this inequality. Then, for every $\varepsilon_2 > 0$, choose $\delta_2(\varepsilon_2) := \delta_1(\varepsilon_2^{1/2}) > 0$. Then,
$$\|L(y^\dagger|\cdot)^{1/2} - L(y^\ddagger|\cdot)^{1/2}\|_{\mathbf{L}^2(X, \muprio)} \leq \varepsilon_2  \quad (y^\ddagger \in Y: \|y^\dagger - y^\ddagger\|_Y \leq \delta_{2}(\varepsilon_2))$$
gives us the desired continuity result.

3. Using the continuity result in 1. and the composition of continuous functions, we also know that $y^\dagger \mapsto Z(y^\dagger)^{-1/2} \in (0,\infty)$ is continuous. Hence, we have for every $y^\dagger \in Y$ and every $\varepsilon_3 > 0$ a $\delta_{3}(\varepsilon_3) > 0$ with 
$$|Z(y^\dagger)^{-1/2} - Z(y^\ddagger)^{-1/2} | \leq \varepsilon_3 \quad (y^\ddagger \in Y: \|y^\dagger- y^\ddagger\|_Y \leq \delta_{3}(\varepsilon_3)).$$
Given this and all the previous results, we now employ a technique that is typically used to prove the continuity of the product of two continuous functions. 
Let $y^\dagger \in Y$, $\varepsilon_2, \varepsilon_3 > 0$, $\delta_4 = \min\{\delta_2(\varepsilon_2), \delta_3(\varepsilon_3)\}$ and $y^\dagger \in Y$: $\|y^\dagger - y^\ddagger\|_Y \leq \delta_4$. 
We arrive at
\begin{align*}
\dHel(\mupost^\dagger,  \mupost^\ddagger) &= \|Z(y^\dagger)^{-1/2}L(y^\dagger|\theta)^{1/2}  - Z(y^\ddagger)^{-1/2}  L(y^\ddagger|\theta)^{1/2}\|_{\mathbf{L}^2(X, \muprio)}\\ 
&\leq |Z(y^\ddagger)^{-1/2}|\times \|L(y^\ddagger|\theta)^{1/2}  - L(y^\dagger|\theta)^{1/2}\|_{\mathbf{L}^2(X, \muprio)} \\ &\quad + \| L(y^\dagger|\theta)^{1/2}\|_{\mathbf{L}^2(X, \muprio)} |Z(y^\ddagger)^{-1/2} - Z(y^\dagger)^{-1/2}|\\
&\leq Z(y^\ddagger)^{-1/2} \varepsilon_{2} + Z(y^\dagger)^{1/2} \varepsilon_3 \\
& \leq ((Z(y^\dagger)^{-1/2} + \varepsilon_3)\varepsilon_2 + Z(y^\dagger)^{1/2} \varepsilon_3 
\end{align*}
where we have used in the last step, that $|Z(y^\dagger)^{-1/2}-Z(y^\ddagger)^{-1/2}| \leq \varepsilon_3$.
We now choose some $\varepsilon_{4} > 0$ and set $\delta_{4} = \min\{\delta_{2}(\varepsilon_2'), \delta_3(\varepsilon_3')\}$, where we set
\begin{align*}
\varepsilon_{2}' := \frac{\varepsilon_4 Z(y^\dagger)^{1/2}}{\varepsilon_4+2}, \quad \varepsilon_3' := \frac{\varepsilon_4}{2Z(y^\dagger)^{1/2}}.
\end{align*}
Then, we obtain that $\dHel(\mupost^\dagger, \mupost^\ddagger) \leq \varepsilon_{4}$ for any $y^\ddagger \in Y$, such that $\|y^\dagger - y^\ddagger\|_Y \leq \delta_{4}$. 
This implies the continuity of the posterior measure in Hellinger distance.
\end{proof}

\begin{proof}[Proof of \Cref{Lem_topologies}]
For every $a \in A$ and $\varepsilon > 0$, there is a $\delta(\varepsilon) > 0$, with $$d_2(f(a), f(a')) \leq \varepsilon \qquad (a' \in A: d_A(a,a') \leq \delta(\varepsilon)).$$
Hence, for the same $a, a', \varepsilon$ and $\delta$, we have 
$$d_1(f(a), f(a')) \leq  t(d_2(f(a), f(a'))) \leq t(\varepsilon)$$
Since $t$ is continuous in $0$, we find for every $\varepsilon' > 0$ some $\delta'(\varepsilon') > 0$, such that $|t(x)| \leq \varepsilon'$ for $x \in [0, \infty): |x| \leq \delta'(\varepsilon')$.
Now, we choose for every $a \in A$ and $\varepsilon'' > 0$: $\delta''(\varepsilon'') := \delta(\delta'(\varepsilon'')).$
Then,
$$d_1(f(a), f(a')) \leq  t(d_2(f(a), f(a'))) \leq t(\delta'(\varepsilon'')) \leq \varepsilon'' \qquad (a' \in A: d_A(a,a') \leq \delta''(\varepsilon''))$$
which results in continuity in $(B, d_1)$.
\end{proof}

\begin{proof}[Proof of \Cref{PropA1A5}]
We show that (A3) and (A5) hold. Note that (A2) is implied by (A3).
We set $g \equiv c$. Then, $L \leq g$. Since, $\muprio$ is a probability measure, we have 
$$\int_X g \mathrm{d}\muprio = c \muprio(X) = c <\infty.$$ 
Hence, $g \in \mathbf{L}^1(X,\muprio)$, which implies that (A3) is satisfied.
Next, we define $g'(\theta') := c \cdot \|\theta'\|_X^p$, for $\theta' \in X, \muprio$-a.s.. By this definition, we have $\|\cdot\|_X^p \cdot L(y^\dagger |\cdot) \leq g'$ for all $y^\dagger \in Y$. Moreover, 
$$
\int_X g' \mathrm{d}\muprio = c \int_X \|\theta\|_X^p \muprio(\mathrm{d}\theta)  <\infty,
$$
since $\muprio \in \Prob_p(X)$ and, thus, $\int_X \|\theta\|_X^p \muprio(\mathrm{d}\theta) <\infty$. Hence, $g' \in \mathbf{L}^1(X,\muprio)$, implying that (A5) holds.
\end{proof}

\begin{proof}[Proof of \Cref{thm:WassWell}]
Let $p \in [1,\infty)$ and $y^\dagger \in Y$. Since (A1)-(A4) hold, we have existence and uniqueness of $\mupost^\dagger \in \Prob(X)$ by \Cref{Thm_main}. We first show, that $\mupost^\dagger \in \Prob_p(X)$:
We have
$$
\int_X \|\theta\|^p_X \mupost^\dagger(\mathrm{d}\theta)= \frac{\int_X L(y^\dagger|\theta) \|\theta\|^p_X \muprio(\mathrm{d}\theta)}{\int_X L(y^\dagger|\theta) \muprio(\mathrm{d}\theta)} \leq \frac{\int_X g'(\theta) \muprio(\mathrm{d}\theta)}{\int_X L(y^\dagger|\theta) \muprio(\mathrm{d}\theta)} < \infty,
$$
where the left-hand side is bounded by \Cref{thm:bayes} (denominator) and by (A5) (numerator). Hence, the posterior measure exists in $\Prob_p(X)$. Since $\Prob_p(X) \subseteq \Prob(X)$, the posterior measure is also unique in $\Prob_p(X)$. Hence, existence and uniqueness of the posterior are satisfied.

Now, we move on to stability.
As in the proof of \Cref{Lemma_main}, the map $Y \ni y^\dagger \mapsto Z(y^\dagger) \in (0, \infty)$ is continuous. By the DCT and (A5), the map $$Y \ni y^\dagger \mapsto \int_X L(y^\dagger|\theta) \|\theta\|^p_X \muprio(\mathrm{d}\theta) \in [0, \infty)$$ is continuous as well.
Therefore,
\begin{align*}
\int_X \|\theta\|^p_X \mupost^\dagger(\mathrm{d}\theta) &= \frac{\int_X L(y^\dagger|\theta) \|\theta\|^p_X \muprio(\mathrm{d}\theta)}{\int_X L(y^\dagger|\theta) \muprio(\mathrm{d}\theta)} \\ &\rightarrow \frac{\int_X L(y^\ddagger|\theta) \|\theta\|^p_X \muprio(\mathrm{d}\theta)}{\int_X L(y^\ddagger|\theta) \muprio(\mathrm{d}\theta)} = \int_X  \|\theta\|^p_X \mupost^\ddagger(\mathrm{d}\theta)
\end{align*}
as $y^\dagger \rightarrow y^\ddagger$.
Hence, we have stability of the posterior measure in the $p$-th moment. 
Additionally, we have weak well-posedness due to \Cref{Thm_main}; and thus stability in the $\dProk$. By \cref{Eq_WasProk}, we have stability in $\dWas$.

Therefore, we also have Wasserstein$(p)$ well-posedness of \cref{EQ_BIP}.
\end{proof}

\begin{proof}[Proof of \Cref{Coro_homeo_wellp}]
According to \Cref{prop_homeomo}, the posterior measure $\mupost^\dagger$ is well-defined and unique.
Let $f:X \rightarrow \mathbb{R}$ be bounded and continuous. Then,
\begin{equation} \label{EQ_Continuous_homeo_proof}
\lim_{y^\ddagger \rightarrow y^\dagger}\int f \mathrm{d}\mupost^\ddagger = \lim_{y^\ddagger \rightarrow y^\dagger}f \circ \mathcal{G}^{-1}(y^\ddagger) = f \circ \mathcal{G}^{-1}(y^\dagger) = \int f \mathrm{d}\mupost^\dagger,
\end{equation}
since $f \circ \mathcal{G}^{-1}$ is continuous. Therefore, $Y \ni y^\dagger \mapsto \mupost^\dagger \in (\Prob(X), \dProk)$ is continuous. Thus, we have weak well-posedness.
If now $X$ is a normed space and $p \in [1,\infty)$, the mapping $\|\cdot\|_X^p: X \rightarrow \mathbb{R}$ is continuous. Note that when setting $f := \|\cdot\|_X^p$ in \Cref{EQ_Continuous_homeo_proof}, the equation still holds. Thus, we have stability in the $p$-th moment and therefore also Wasserstein$(p)$ well-posedness according to \cref{Eq_WasProk}.
\end{proof}

\begin{proof}[Proof of \Cref{thm:KLD}]
First note that (A1)-(A4) imply the existence and the uniqueness of the posterior measure, as well as the continuity of $y^\dagger \mapsto Z(y^\dagger)$.
Let $y^\dagger \in Y$ and  $y^\ddagger \in Y$, with $\|y^\dagger - y^\ddagger \|_Y \leq \delta$. $\delta > 0$ is chosen as in (A9). We have
\begin{align*}
\DKL(\mupost^\dagger\|\mupost^\ddagger) &= \int \log\left(\frac{\mathrm{d}\mupost^\dagger}{\mathrm{d}\mupost^\ddagger}\right) \mathrm{d}\mupost^\dagger \\ &= \int  \log L(y^\dagger|\cdot)  -\log L(y^\ddagger|\cdot)\mathrm{d}\mupost^\dagger +\left(\log Z(y^\ddagger) - \log Z(y^\dagger) \right),
\end{align*}
where the right-hand side of this equation is well-defined since $Z(y^\dagger), Z(y^\ddagger) \in (0, \infty)$ by \Cref{Lemma_main} and since (A9) holds. 
Moreover, the continuity in the model evidence implies that $\left(\log Z(y^\ddagger) - \log Z(y^\dagger) \right) \rightarrow 0$, as $y^\ddagger \rightarrow y^\dagger$. 
Also, note that $\log L(\cdot \, |\theta')$ is continuous by (A4), which implies
\begin{align*}
\lim_{y^\ddagger \rightarrow y^\dagger} \int_X \log L(y^\dagger|\cdot)  -\log L(y^\ddagger|\cdot) \mathrm{d}\mupost^\dagger = \int_X \lim_{y^\ddagger \rightarrow y^\dagger} \log L(y^\dagger|\cdot)  -\log L(y^\ddagger|\cdot) \mathrm{d}\mupost^\dagger = 0,
\end{align*}
where we applied the DCT with $2h(\cdot ,y^\dagger)$ as a dominating function.
\end{proof}

\begin{proof}[Proof of \Cref{Cor_Gaussian_case}]
We check (A1)-(A4). (A1) By definition, the likelihood is a strictly positive probability density function for any $\theta' \in X$.  (A2)-(A3) The likelihood is bounded above uniformly by $g \equiv \det(2\pi\Gamma)^{-1/2}$ which is integrable with respect to any probability measure on $(X,\mathcal{B}X)$. (A4) the likelihood is  continuous in $y^\dagger$ for any $\theta' \in X$.
\end{proof}

\begin{proof}[Proof of \Cref{Prop_infinite_dim_Gauss}]
1. The function $L$ is indeed a correct likelihood, i.e.  $y\mapsto L(y^\dagger|\theta')$ is probability density function for $\muprio$-a.e. $\theta' \in X$. We refer to the discussions of the Cameron--Martin Theorem in \cite[\S 2.4]{Bogachev1998} and \cite[\S 2.7]{Sullivan2015}. Moreover, we again mention \cite[Remark 3.8]{Stuart2010} and  \cite[\S 2.1]{Kahle2019} who have discussed the modelling in this case. Hence, (A1) is true.

2. Now, we check (A2)-(A4). Note that (A4) is true by assumption. (A2) holds since $\mathcal{G}$ is bounded. (A3) cannot be shown easily. However, we can replace it by a local version of this assumption; see \Cref{Remark_local_versions_A}. Indeed, to show continuity of $\mupost^\dagger$ in $y^\dagger \in Y$, we only need to satisfy (A3) in $\overline{B}(y^\dagger, \delta) := \{\|\cdot - y^\dagger\|_Y \leq \delta\} := \{y^\ddagger : \|y^\ddagger -y^\dagger \|_Y \leq \delta \},$
for $\delta > 0$. If we show this for any $y^\dagger \in Y$ and some $\delta > 0$, we obtain stability as well. 
Note that we have used this idea to show Kullback--Leibler well-posedness in \Cref{thm:KLD}.

3. Let $y^\dagger \in Y$ be arbirary. Let $c$ be chosen such that $\|\mathcal{G}(\theta')\|_Y < c$, which exists since $\mathcal{G}$ is bounded. Let $y^\ddagger \in \overline{B}(y^\dagger, \delta)$. By the Cauchy--Schwarz and the triangle inequality, we have
\begin{align*}
L(y^\ddagger|\theta') &= \exp\left(\langle \mathcal{G}(\theta'), y^\ddagger\rangle_Y - \frac{1}{2}\|\mathcal{G}(\theta')\|_Y^2 \right) \leq \exp\left( |\langle \mathcal{G}(\theta'), y^\ddagger\rangle_Y|  \right) 
\\ &\leq \exp\left(\| \mathcal{G}(\theta')\|_Y \|y^\ddagger \|_Y \right) \leq \exp(c \|y^\ddagger\|_Y) \\  &= \exp(c \|y^\ddagger - y^\dagger + y^\dagger\|_Y) \leq  \exp(c \|y^\ddagger - y^\dagger\|_Y + \|y^\dagger\|_Y) \\ &\leq  \exp(c \cdot (\delta + \| y^\dagger\|_Y)) =: c',
\end{align*}
for $\muprio$-almost every $\theta' \in X$.
Now, we choose $g: X \rightarrow \mathbb{R}$  to be $g \equiv c'$.
Let now $\muprio \in \mathrm{Prob}(X)$.
 Then, $g \in \mathbf{L}^1(X,\muprio)$ and 
$L(y^\ddagger| \theta') \leq g(\theta')$, for $\muprio$-almost every $\theta' \in X$. Since $y^\dagger$ is chosen arbitrarily, we obtain stability in the weak topology, the Hellinger distance, and the total variation distance.
Hence, we have whtv well-posedness and thus we have shown (a).

4. Let $p \in [1,\infty)$. To show Wasserstein$(p)$ well-posedness, we can again use a local argument on the data space. Hence, we can satisfy (A5) locally on the data space. This on the other hand is implied by a local version of \Cref{PropA1A5}. Hence, we obtain stability in the Wasserstein$(p)$ distance, if $\muprio \in \Prob_p(X)$ and if for all $y^\dagger \in Y$, we have some $\delta > 0$ and $c' > 0$ such that
$$
L(y^\ddagger|\theta') \leq c' \qquad (\|y^\ddagger - y^\dagger \|_Y \leq \delta; \muprio\text{-a.e. } \theta' \in X).
$$
This however is what we have shown already in 3. Thus, we have shown (b).
\end{proof}

\begin{proof}[Proof of \Cref{Thm_RegCondProb}]
The theorem above holds, if $\Omega, X, Y$ are Radon spaces; see \cite[Theorem 3.1]{Leao2004}. $Y$ is a Radon space by definition. $X$ and $\Omega$ can be extended to Radon spaces $X'$ and $\Omega' = X' \times Y$, where $\mathbb{P}(\theta \in X' \backslash X) = 0 = \mathbb{P}(\Omega' \backslash \Omega)$. Moreover, we set $M(y,X' \backslash X) = 0$ $(y \in Y)$.
\end{proof}

\begin{proof}[Proof of \Cref{Lemma_marginalpdf}] Let $A \in \mathcal{B}X$.
Note that $$\mathbb{P}(\theta \in A) = \mathbb{P}((\theta, y) \in A \times Y) = \int_{A \times Y} f \mathrm{d}(\nu_X \otimes \nu_Y) = \int_A \int_Y f(\theta^\dagger, y^\dagger)\nu_Y(\mathrm{d}y^\dagger)\nu_X(\mathrm{d}\theta^\dagger),$$
where the last equality holds due to Tonelli.
Hence, indeed 
$$\frac{\mathrm{d}\mathbb{P}(\theta \in \cdot)}{\mathrm{d}\nu_X} = \int_X f(\cdot, y^\dagger) \nu_Y(\mathrm{d}y^\dagger) \quad (\nu_X\text{-a.e.}).$$
The statement about the $\nu_Y$-probability density function of $y$ can be shown by exchanging $y$ and $\theta$, and $X$ and $Y$.
\end{proof}

\begin{proof}[Proof of \Cref{Lemma_condpdf}]
For a derivation in the case $X := Y := \mathbb{R}$, see \cite[Example 5.3.2 (b)]{Ash2000}. The proof in our more general setting is analoguous.
\end{proof}
\section*{Acknowledgments}
The author thanks several contributors for their highly appreciated support:
Elisabeth Ullmann made insightful and valuable comments that contributed to this work, as did illuminating discussions with Bj\"orn Sprungk. Tim J. Sullivan detected an error in an older version of this manuscript.
Florian Beiser and Brendan Keith have proofread this article. 
Lukas Latz helped the author interpreting Hadamard's original work.
Both anonymous reviewers made valuable comments which helped to improve this manuscript.
\bibliographystyle{siamplain}
\bibliography{library}

\end{document}